\documentclass[10pt]{article}
\usepackage{amsmath}
\usepackage{amsfonts}
\usepackage{amssymb}
\usepackage{paralist}
\usepackage{graphics} %% add this and next lines if pictures should be in esp format
\usepackage{epsfig} %For pictures: screened artwork should be set up with an 85 or 100 line screen
\usepackage[colorlinks=true]{hyperref}

\hypersetup{urlcolor=blue, citecolor=red}

\newtheorem{theorem}{Theorem}[section]
\newtheorem{corollary}{Corollary}[section]

\newtheorem{lemma}{Lemma}[section]
\newtheorem{proposition}{Proposition}[section]

\newtheorem{definition}{Definition}[section]
\newtheorem{remark}{Remark}[section]

\newcommand{\ep}{\varepsilon}

\def\sqr#1#2{{\vcenter{\vbox{\hrule height.#2pt
              \hbox{\vrule width.#2pt height#1pt \kern#1pt \vrule
width.#2pt}
              \hrule height.#2pt}}}}
\def\3n{\negthinspace \negthinspace \negthinspace }
\def\2n{\negthinspace \negthinspace }
\def\1n{\negthinspace }

\def\ds{\displaystyle}
\def\cF{{\mathcal F}}

\def\q{\quad}
\def\qq{\qquad}

\def\({\Big (}
\def\){\Big )}
\def\[{\Big[}
\def\]{\Big]}

\makeatletter
   
   \@addtoreset{equation}{section}
\makeatother

\def\ba{\begin{array}}
\def\be{\begin{equation}}
\def\bc{\begin{corollary}}
\def\bd{\begin{definition}}
\def\bea{\begin{eqnarray}}
\def\bel{\begin{equation}\label}
\def\bl{\begin{lemma}}
\def\bp{\begin{proposition}}
\def\br{\begin{remark}}
\def\bt{\begin{theorem}}
\def\ee{\end{equation}}
\def\eea{\end{eqnarray}}
\def\et{\end{theorem}}
\def\ec{\end{corollary}}
\def\el{\end{lemma}}
\def\ep{\end{proposition}}
\def\er{\end{remark}}
\def\ea{\end{array}}
\def\ed{\end{definition}}

\title%[LQ Time-Inconsistent Optimal Control]
      {Discrete time mean-field stochastic linear-quadric optimal control problems\thanks{Professor Robert Elliott wishes to thank the ARC and NSERC for support.
Dr. Xun Li wishes to thank the Hong Kong RGC for support with grant no. 520412.
Dr. Yuan-Hua Ni wishes to thank the support from the National Natural Science Foundation of China under grant 11101303. Corresponding author: Yuan-Hua Ni. Tel. +86-22-27427090.}}

\author{Robert  Elliott\thanks{School of Mathematical Sciences, University of Adelaide, SA 5005, Australia; Haskayne School of Business, University of Calgary, Calgary, Alberta, Canada ({\tt relliott@ucalgary.ca})}~~~~~~Xun Li\thanks{Department of Applied Mathematics, The Hong Kong Polytechnic University, Hunghom, Kowloon, Hong Kong ({\tt malixun@inet.polyu.edu.hk})}~~~~~~Yuan-Hua Ni\thanks{Department of Mathematics, School of Science, Tianjin Polytechnic University, Tianjin, P.R. China ({\tt yhni@amss.ac.cn}). }}

\begin{document}
\maketitle

%% Enter the first author's name and address:
%\centerline{\scshape Jiongmin Yong }
%\medskip
%{\footnotesize
%% please put the address of the first author
% \centerline{Department of Mathematics, University of Central Florida}
%   \centerline{ Orlando, FL 32816, USA}
%} % Do not forget to end the {\footnotesize by the sign }
%
%
%\bigskip
%
%% The name of the associate editor will be entered by an editorial staff
%% "Communicated by the associate editor name" is not needed for special issue.
% \centerline{(Communicated by Gengsheng Wang)}

%The abstract of your paper
\begin{abstract}

This paper first presents necessary and sufficient conditions for the solvability of discrete time, mean-field, stochastic linear-quadratic optimal control problems. Then, by introducing several sequences of bounded linear operators, the  problem becomes an operator stochastic LQ problem, in which the optimal control is a linear state feedback. Furthermore, from the form of the optimal control, the problem changes to a matrix dynamic optimization problem. Solving this optimization problem, we obtain the optimal feedback gain and thus the optimal control. Finally, by completing the square, the optimality of the above control is validated.
\end{abstract}

%The title of your section 1
\section{Introduction}{\label 1}

In this paper we consider a class of stochastic linear-quadratic (LQ) optimal control problems of mean-field type. The system equation is the following linear stochastic difference equation with $k\in\{0,1,2,...,N-1\}\equiv \mathbb{N}$,
\begin{eqnarray}\label{system}
\left\{\begin{array}{ll}
x_{k+1}=(A_kx_k+\bar{A}_k\mathbb{E}x_k+B_ku_k+\bar{B}_k\mathbb{E}u_k)+(C_kx_k+\bar{C}_k\mathbb{E}x_k+D_ku_k+\bar{D}_k\mathbb{E}u_k)w_k, \\
x_0=\zeta,
\end{array}\right.
\end{eqnarray}
where $A_k,\bar{A}_k,C_k,\bar{C}_k\in \mathbb{R}^{n\times n}$, and $B_k,\bar{B}_k,D,\bar{D}_k\in \mathbb{R}^{n\times m}$ are given deterministic matrices, and $\mathbb{E}$ is the expectation operator. Denote  $\{0,1,2,...,$ $ N\}$ by $\bar{\mathbb{N}}$. In (\ref{system}) , $\{x_k, k\in \bar{\mathbb{N}}\}$ and $\{u_k, k\in \mathbb{N}\}$ are the state process and control process, respectively. $\{w_k, k\in \mathbb{N}\}$, defined on a probability space $(\Omega, \mathcal{F}, P)$, represents the stochastic disturbances, which is assumed to be a martingale difference sequence
\begin{eqnarray}\label{w-moment-0}
\mathbb{E}[w_{k+1}|\cF_k]=0, \quad
\mathbb{E}[(w_{k+1})^2|\cF_k]=1,
\end{eqnarray}
where $\mathcal{F}_k$ is the $\sigma$-algebra generated by $\{\zeta,
w_l, l=0, 1,\cdots,k\}$. The initial value $\zeta$ and $\{w_k, k\in \mathbb{N}\}$ are assumed to be independent of each other. The cost functional associated with (\ref{system}) is
\begin{eqnarray}\label{cost-finite}
\begin{array}{rl}
J(\zeta,u) =  & \ds\mathbb{E}\left[\sum_{k=0}^{N-1}\left(x_{k}^TQ_kx_{k}+(\mathbb{E}x_{k})^T\bar{Q}_k\mathbb{E}x_{k}+u_k^TR_ku_k+(\mathbb{E}u_{k})^T\bar{R}_k\mathbb{E}u_{k}\right) \right] \\ [5mm]
&+\mathbb{E}\left(x_{N}^T{G}_Nx_{N}\right)+(\mathbb{E}x_{N})^T\bar{G}_N\mathbb{E}x_{N},
\end{array}
\end{eqnarray}
where $Q_k, \bar{Q}_k, R_k, \bar{R}_k, k\in\mathbb{N}, G_N, \bar{G}_N$ are deterministic symmetric matrices with appropriate dimensions.

We introduce the following admissible control set
\begin{eqnarray*}\label{admissilbe-control-set-finite}
\begin{array}{l}
\mathcal{U}_{ad}\equiv\left\{u=(u_0,u_1,\cdots,u_{N-1}):~\mathbb{N}\times\Omega\rightarrow \mathbb{R}^m~\left\vert~u_k\in \mathcal{F}_k,\mathbb{E}|u_k|^2<\infty\right.\right\}.
\end{array}
\end{eqnarray*}
The optimal control problem considered in
this paper is stated as follows:

\noindent
\textbf{Problem (MF-LQ)}.
\emph{For any given square-integrable initial value $\zeta$, find  $u^o\in \mathcal{U}_{ad}$ such that}
\begin{eqnarray} \label{optimal-finite}
J(\zeta,u^o)=\inf_{u\in \mathcal{U}_{ad}}J(\zeta,u).
\end{eqnarray}
\emph{We then call $u^o$ an optimal control for Problem (MF-LQ).}

Unlike the classical stochastic LQ problem, the expectation $\mathbb{E}x_k$ of $x_k$ appears in the system equation (\ref{system}) and the cost functional (\ref{cost-finite}). Therefore, the above problem is called a mean-field LQ problem. This problem is a combination of mean field theory and a LQ problem. Mean-field theory was developed to study the collective behaviors resulting from individuals' mutual interactions in various physical and sociological dynamic systems.
According to mean-field theory, the interactions among agents are modelled by a mean-field term. Letting the number of individuals go to the infinity, the mean-field term will approach the expectation. To see this, assume the dynamics of particle $i$ $(i=1,\cdots,L)$ are
\begin{eqnarray*}
x^{i,L}_{k+1}=\Big{(}A_kx^{i,L}_k+\bar{A}_k\frac{1}{L}\sum_{j=1}^Lx^{j,L}_k+B_ku_k\Big{)}+\Big{(}C_kx^{i,L}_k+\bar{C}_k\frac{1}{L}\sum_{j=1}^Lx^{j,L}_k+D_ku_k\Big{)}w_k.
\end{eqnarray*}
Under appropriate conditions and letting $L \rightarrow \infty$, we have
\begin{eqnarray}\label{system-limit}
x_{k+1}=(A_kx_k+\bar{A}_k\mathbb{E}x_k+B_ku_k)+(C_kx_k+\bar{C}_k\mathbb{E}x_k+D_ku_k)w_k.
\end{eqnarray}
An exact derivation of (\ref{system-limit}) follows the classical
Mckean-Vlasov argument; for the Mckean-Vlasov argument, readers may
see, for example, \cite{Bensoussan 2012}\cite{Mckean}\cite{Sznitman}
and references therein. Clearly, (\ref{system}) is a natural
extension of (\ref{system-limit}). For the motivation for including
$\mathbb{E}x_k$ and $\mathbb{E}u_k$ in the cost functional
(\ref{cost-finite}), \cite{Yong-MFF-2011} points out that it is
natural to introduce variations $\mbox{var}(x_k)$ and
$\mbox{var}(u_k)$ to the cost functional so as to the state process
and the control process could be not too sensitive to the random
events. An example of this case is the known finance mean-variance
problem. For system (\ref{system}) without mean-field terms, we
refer to papers such as
\cite{Ait-Chen-Zhou-2002}\cite{Zhang-weihai-1}\cite{Morozan}\cite{Zhang-weihai-2}.

The continuous-time counterpart of (\ref{system}) is a mean-field
stochastic differential equation (SDE), whose investigation goes
back to the McKean-Vlasov SDE proposed in 1950s
(\cite{Kac}\cite{Mckean}). Since then,
many researches study
McKean-Vlasov type SDEs and applications (\cite{Chan}\cite{Dawson}\cite{Gartner}\cite{Graham}). For recent development of
mean-field SDEs, readers may refer to
\cite{Buckdahn-1}\cite{Buckdahn-2}\cite{Crisan} and
 references therein. The control problems for
mean-field SDEs are investigated by many authors; see, for example,
\cite{Ahmed-2}\cite{Ahmed}\cite{Andersson}\cite{Bensoussan
2012}\cite{Buckdahn-3}\cite{zhou}\cite{Yong-MFF-2011}\cite{Yong-MFI-2012}.
On the other hand, another class of problems is
mean-field games in terms of its ability to model the collective
behavior of individuals due to their mutual interactions. Comparing
the above mentioned papers, this class of problems may be viewed as
decentralized control problems, that is, the controls are selected
to achieve each individual's own goal. For other aspects of mean
field games, readers may refer to \cite{Bensoussan
2012}\cite{Huang2003}\cite{Huang2007}\cite{Huang2012}\cite{Lasry}\cite{Li-tao}\cite{Li-tao2}\cite{Wang
bingchang} and references therein.

The problem considered in this paper is related to that of
\cite{Yong-MFF-2011}, which deals with a continuous time mean-field
LQ problem. In \cite{Yong-MFF-2011}, by a variational method,
the optimality system is derived, which is a mean-field
forward-backward SDE.
Furthermore, by a decoupling technique, two Riccati differential
equations are obtained.  This gives the feedback representation of
the optimal control. Our discussion of the mean-field LQ optimal
control problem differs from \cite{Yong-MFF-2011} in the following
ways. First, we considers the discrete time case. There are several
situations when it is necessary, or natural, to describe a system by
a discrete time model. A typical case is when the signal values are
only available for measurement or manipulation at certain times, for
example, because a continuous time system is sampled at certain
times. Another case arises from discretization of the dynamics of
continuous time problems. The second difference between this paper
and \cite{Yong-MFF-2011} is that methodology differs. In this paper,
we first convert Problem (MF-LQ) to a quadratic optimization problem
in Hilbert space. This gives necessary and sufficient conditions for
the solvability of Problem (MF-LQ). The optimal control obtained in
this case has an abstract form, which may be viewed as
open loop. Secondly, by introducing several sequences
of bounded linear operators, Problem (MF-LQ) becomes an operator LQ
optimal control problem. This operator LQ problem shows that the
optimal control is a linear state feedback, which is clearly
closed-loop. Thirdly, by the linearity of the optimal control,
Problem (MF-LQ) becomes a matrix dynamic optimization problem. Using
the matrix minimum principle (\cite{Athans}), we derive the optimal
feedback gains and thus the optimal control is obtained. Finally, by
completing the square, we validate the correctness of the obtained
optimal control.

The paper is organized as follows. The next section presents the preliminaries and two abstract considerations---quadratic optimization in Hilbert space and the operator LQ problem. In Section 3, optimal control via Riccati equations is presented. Section 4 gives an example to calculate the solutions to Riccati equations and the optimal feedback gains. Section 5 gives some conclude remarks.

\section{Preliminaries and abstract considerations}

Some standard notion is introduced.
\begin{definition}
(i). Problem (MF-LQ) is said to be finite for $\zeta$ if
\begin{eqnarray*}
\inf_{u\in \mathcal{U}_{ad}}J(\zeta,u)>-\infty.
\end{eqnarray*}
Problem (MF-LQ) is said to be finite if it is finite for any $\zeta$.

(ii). Problem (MF-LQ) is said to be (uniquely) solvable for $\zeta$ if there exists a (unique) $u^o\in \mathcal{U}_{ad}$
such that (\ref{optimal-finite}) holds for $\zeta$. Problem (MF-LQ) is said to be (uniquely) solvable if it is solvable for any $\zeta$.
\end{definition}

In this subsection, we shall consider Problem (MF-LQ) using two methods. The first converts Problem (MF-LQ) to a quadratic optimization problem in Hilbert space, which gives necessary and sufficient conditions for the solvability of Problem (MF-LQ). The second considers Problem (MF-LQ) in the language of an operator LQ problem. This reveals that the optimal control is a linear state feedback.

Introduce the following spaces:
\begin{eqnarray*}
\begin{array}{l}
\mathcal{X}_k=L_{\mathcal{F}_{k}}^{2}(\mathbb{R}^n)=\left\{\xi:\Omega\mapsto \mathbb{R}^n| \xi\mbox{ is }\mathcal{F}_{k}\mbox{-measurable}, \mathbb{E}|\xi|^2<\infty \right\},~~k\in\bar{\mathbb{N}},\\ [5mm]
\mathcal{X}[0,k]=\left\{(x_0,\cdots,x_k) \left\vert x_k\in\mathcal{X}_k, \mbox{ and is }\mathcal{F}_k\mbox{-measurable}, \sum_{l=0}^{k}\mathbb{E}|x_l|^2<\infty \right.\right\},~~k\in\bar{\mathbb{N}}, \\ [5mm]
\mathcal{U}_k=L_{\mathcal{F}_{k}}^{2}(\mathbb{R}^m)=\left\{\eta:\Omega\mapsto \mathbb{R}^m| \eta\mbox{ is }\mathcal{F}_{k}\mbox{-measurable}, \mathbb{E}|\eta|^2<\infty \right\}, ~~k\in \mathbb{N}.
\end{array}
\end{eqnarray*}
Clearly, for any $k\in \bar{\mathbb{N}}, l\in \mathbb{N}$, $ \mathcal{X}[0,k]$ and $\mathcal{H}=\mathcal{X}_k, \mathcal{U}_l$ are Hilbert spaces under the usual inner products
\begin{eqnarray}\label{inner-product}
\langle x, z\rangle=\mathbb{E}\left(x^Tz \right),\mbox{  for any }x,z\in\mathcal{H},
\end{eqnarray}
and
\begin{eqnarray}
\langle x, z\rangle=\mathbb{E}\left(\sum_{p=0}^kx_p^Tz_p \right),~~~\mbox{  for any }x=(x_0,...,x_k),z=(z_0,...,z_k)\in\mathcal{X}[0,k].
\end{eqnarray}
For any variable $z$ in $\mathcal{X}_k$ or $\mathcal{U}_k$, the expectation of $z$, i.e., $\mathbb{E}z$, is clearly well defined. If we consider $\mathbb{E}$ as an operator, it is clear that the domain and range of $\mathbb{E}$ may differ from place to place. For example, the domain may be $\mathcal{X}_k$, $\mathcal{U}_k$, and the range is $\mathbb{R}^n$ and $\mathbb{R}^m$, respectively. Therefore, the domain and range of the adjoint operator $\mathbb{E}^*$ of $\mathbb{E}$ may differ. For example, if $\mathbb{E}$ maps $\mathcal{X}_k$ to $\mathbb{R}^n$, then $\mathbb{E}^*$ is defined from $\mathbb{R}^n$ to $\mathcal{X}_k$ and is defined as
\begin{eqnarray*}
\langle \mathbb{E}^*r_n, x\rangle=\langle r_n, \mathbb{E}x\rangle,~~~~\mbox{for any }r_n\in\mathbb{R}^n, x\in \mathcal{X}_k.
\end{eqnarray*}
If $\mathbb{E}$ maps $\mathcal{U}_k$ to $\mathbb{R}^m$, then $\mathbb{E}^*$ is defined from $\mathbb{R}^m$ to $\mathcal{U}_k$, and is defined as
\begin{eqnarray*}
\langle \mathbb{E}^*r_m, z\rangle=\langle r_m, \mathbb{E}z\rangle,~~~~\mbox{for any }r_m\in\mathbb{R}^m, z\in \mathcal{U}_k.
\end{eqnarray*}
Consequently, $\mathbb{E}$ and $\mathbb{E}^*$ could be denoted by $\mathbb{E}_\mathcal{H}$ and $\mathbb{E}^*_\mathcal{H}$ in order to emphasize $\mathcal{H}$, which is the domain of $\mathbb{E}$ and thus the range of $\mathbb{E}^*$. In this paper, $\mathcal{H}$ may be $\mathcal{X}_k, \mathcal{U}_k, k\in \mathbb{N}$, and  $\mathbb{E}^*_\mathcal{H}$ always appears accompanying $\mathbb{E}_\mathcal{H}$ in the form of $\mathbb{E}^*_\mathcal{H}K\mathbb{E}_\mathcal{H}$ with $K$ being a generic square matrix. Clearly, for any $z\in \mathcal{H}$, $\mathbb{E}_\mathcal{H}^*K\mathbb{E}_\mathcal{H}z$ is in $\mathcal{H}$. Another form is $\mathbb{E}^*_{\mathcal{U}_k}M\mathbb{E}_{\mathcal{X}_k}$, $k\in\mathbb{N}$, with $M\in \mathbb{R}^{m\times n}$. To simplify the notation, throughout this paper, $\mathbb{E}_\mathcal{H}^*K\mathbb{E}_\mathcal{H}z$ with $z\in\mathcal{H}$, and $\mathbb{E}^*_{\mathcal{U}_k}M\mathbb{E}_{\mathcal{X}_k}x$ with $x\in\mathcal{X}_k$, will be denoted by $\mathbb{E}^*K\mathbb{E}z$ and $\mathbb{E}^*M\mathbb{E}x$, respectively. The meanings will be understood during the context.

Taking expectations in (\ref{system}), we have
\begin{eqnarray}\label{system-Ex}
\left\{\begin{array}{l}
\mathbb{E}x_{k+1}=(A_k+\bar{A}_k)\mathbb{E}x_k+(B_k+\bar{B}_k)\mathbb{E}u_k,~~k\in \mathbb{N}, \\ [2mm]
\mathbb{E}x_0=\mathbb{E}\zeta.
\end{array}\right.
\end{eqnarray}
Let
\begin{eqnarray*}
\left\{\begin{array}{l}
\bar{\Phi}(k,l)=(A_{k}+\bar{A}_{k})(A_{k-1}+\bar{A}_{k-1})\cdots (A_{l}+\bar{A}_{l}),~~k\geq l, \\ [2mm]
\bar{\Phi}(k,l)=I,~~k<l.
\end{array}\right.
\end{eqnarray*}
Then we have
\begin{eqnarray*}
\mathbb{E}x_{k+1}=\bar{\Phi}(k,0)\mathbb{E}\zeta+\sum_{l=1}^{k}\bar{\Phi}(k, l)(B_{l-1}+\bar{B}_{l-1})\mathbb{E}u_{l-1},~~k\in\mathbb{N}\setminus\{0\}.
\end{eqnarray*}
On the other hand, let
\begin{eqnarray*}
\left\{\begin{array}{l}
{\Phi}(k,l)=(A_{k}+w_{k}C_{k})(A_{k-1}+w_{k-1}C_{k-1})\cdots (A_{l-1}+w_{l-1}C_{l-1}),~~ k\geq l, \\ [2mm]
{\Phi}(k,l)=I, k<l.
\end{array}
\right.\end{eqnarray*}
By (\ref{system}), we have for $k\in \mathbb{N}$
\begin{eqnarray*}
\begin{array}{rl}
x_{k+1} = & {\Phi}(k,0)\zeta +\ds\sum_{l=1}^{k}\Phi(k,l)\left[(\bar{A}_{l-1}+w_{l-1}\bar{C}_{l-1})\mathbb{E}x_{l-1}+(B_{l-1}+w_{l-1}D_{l-1})u_{l-1}+(\bar{B}_{l-1}+w_{l-1}\bar{D}_{l-1})\mathbb{E}u_{l-1}\right] \\ [2mm]
= & \ds{\Phi}(k,0)\zeta+\sum_{l=1}^{k}\Phi(k,l)(\bar{A}_{l-1}+w_{l-1}\bar{C}_{l-1})\bar{\Phi}(l-2,0)\mathbb{E}\zeta \\ [2mm]
&+\ds\sum_{l=1}^{k}\left[\Phi(k,l)(\bar{A}_{l-1}+w_{l-1}\bar{C}_{l-1})\sum_{i=1}^{l-2}\bar{\Phi}(l-2, i)(B_{i-1}+\bar{B}_{i-1})\mathbb{E}u_i\right] \\ [2mm]
&+\ds\sum_{l=1}^{k}\Phi(k,l)(B_{l-1}+w_{l-1}D_{l-1})u_{l-1}+\sum_{l=1}^{k}\Phi(k,l)(\bar{B}_{l-1}+w_{l-1}\bar{D}_{l-1})\mathbb{E}u_{l-1}.
\end{array}
\end{eqnarray*}
Now define the following operators for any $\zeta\in \mathcal{X}_0, u\in \mathcal{U}_{ad}$:
\begin{eqnarray*}
\left\{\begin{array}{l}
(\Gamma \zeta)(\cdot)={\Phi}(\cdot-1,0)\zeta, \\ [2mm]
\hat{\Gamma} \zeta=(\Gamma \xi)(N), \\ [2mm]
(\bar{\Gamma}\zeta)(\cdot)=\sum_{l=1}^{\cdot-1}\Phi(\cdot-1,l)(\bar{A}_{l-1}+w_{l-1}\bar{C}_{l-1})\bar{\Phi}(l-2,0)\mathbb{E}\zeta, \\ [2mm]
\hat{\bar{\Gamma}}\zeta=(\bar{\Gamma}\zeta)(N), \\ [2mm]
(Lu)(\cdot)=\sum_{l=1}^{\cdot-1}\Phi(\cdot-1,l)(B_{l-1}+w_{l-1}D_{l-1})u_{l-1}, \\ [2mm]
\hat{L}u= (Lu)(N), \\ [2mm]
(\bar{L}u)(\cdot)=\sum_{l=1}^{\cdot-1}\Phi(\cdot-1,l)(\bar{A}_{l-1}+w_{l-1}\bar{C}_{l-1})\sum_{i=1}^{l-2}\bar{\Phi}(l-2, i)(B_{i-1}+\bar{B}_{i-1})\mathbb{E}u_{i-1} \\ [2mm]
\qq\qq +\sum_{l=1}^{\cdot-1}\Phi(\cdot-1,l)(\bar{B}_{l-1}+w_{l-1}\bar{D}_{l-1})\mathbb{E}u_{l-1},\\
\hat{\bar{L}}u=(\bar{L}u)(N).
\end{array}\right.
\end{eqnarray*}
Then
\begin{eqnarray*}
x_k=(\Gamma \zeta)(k)+(\bar{\Gamma} \zeta)(k)+(Lu)(k)+(\bar{L}u)(k).
\end{eqnarray*}
Clearly, the operators
\begin{eqnarray}\label{operators}
\left\{\begin{array}{l}
\Gamma, \bar{\Gamma}: \mathcal{X}_{0}\mapsto \mathcal{X}[0,N],~~\hat{\Gamma}: \mathcal{X}_{0}\mapsto \mathcal{X}_N,\\
L, \bar{L}: \mathcal{U}_{ad}\mapsto \mathcal{X}[0,N],~~\hat{\bar{L}}: \mathcal{U}_{ad}\mapsto \mathcal{X}_{N},
\end{array}\right.
\end{eqnarray}
are all bounded and linear. Notice that the spaces in (\ref{operators}) are all Hilbert space. Therefore, the corresponding adjoint operators uniquely exist.
Further, in what follows, we use the convention
\begin{eqnarray*}
\left\{\begin{array}{l}
(Qx)(\cdot)=Q_kx_\cdot,~~\forall x\in \mathcal{X}[0,N-1],\\
(\bar{Q}\varphi)(\cdot)=\bar{Q}_k\varphi_\cdot,~~\forall \varphi=(\varphi_0,\cdots,\varphi_{N-1})\mbox{ with }\varphi_k\in\mathbb{R}^n \mbox{ such that }\sum_{k=0}^{N-1}|\varphi_k|^2<\infty,\\
(Ru)(\cdot)=R_ku_\cdot,~~\forall u\in \mathcal{U}_{ad},\\
(\bar{R}\psi)(\cdot)=\bar{R}_k\psi_\cdot,~~\forall
\psi=(\psi_0,\cdots,\psi_{N-1})\mbox{ with }\psi_k\in\mathbb{R}^m
\mbox{ such that }\sum_{k=0}^{N-1}|\psi_k|^2<\infty.
\end{array}\right.
\end{eqnarray*}

Consequently, the cost functional $J(\zeta,u)$ has the following form
\begin{eqnarray*}
\begin{array}{rl}
J(\zeta,u)= & \langle Q(\Gamma \zeta+\bar{\Gamma}\zeta+Lu+\bar{L}u), \Gamma \zeta+\bar{\Gamma}\zeta+Lu+\bar{L}u\rangle+\langle Ru, u\rangle \\ [2mm]
& +\langle \bar{Q}\mathbb{E}(\Gamma \zeta+\bar{\Gamma}\zeta+Lu+\bar{L}u), \mathbb{E}(\Gamma \zeta+\bar{\Gamma}\zeta+Lu+\bar{L}u)\rangle+\langle \bar{R}\mathbb{E}u, \mathbb{E}u \rangle \\ [2mm]
& +\langle G_T(\hat{\Gamma} \zeta+\hat{\bar{\Gamma}} \zeta+\hat{L}u+\hat{\bar{L}}u), \hat{\Gamma} \zeta \\ [2mm]
& +\hat{\bar{\Gamma}} \zeta+\hat{L}u+\hat{\bar{L}}u\rangle+\langle \bar{G}_T\mathbb{E}(\hat{\Gamma} \zeta+\hat{\bar{\Gamma}} \zeta+\hat{L}u+\hat{\bar{L}}u), \mathbb{E}(\hat{\Gamma} \zeta+\hat{\bar{\Gamma}} \zeta+\hat{L}u+\hat{\bar{L}}u)\rangle \\ [2mm]
= & \langle \Theta_1u, u\rangle+2\langle \Theta_2\xi, u\rangle+\langle \Theta_3 \xi, \xi\rangle.
\end{array}
\end{eqnarray*}
Here
\begin{eqnarray*}
\begin{array}{l}
\Theta_1=R+\mathbb{E}^*\bar{R}\mathbb{E}+(L+\bar{L})^*Q(L+\bar{L})+(L+\bar{L})^*\mathbb{E}^*Q\mathbb{E}(L+\bar{L}) +(\hat{L}+\hat{\bar{L}})^*G(\hat{L}+\hat{\bar{L}}) \\ [2mm]
\qq +(\hat{L}+\hat{\bar{L}})^*\mathbb{E}^*G\mathbb{E}(\hat{L}+\hat{\bar{L}}), \\ [2mm]
\Theta_2=(L+\bar{L})^*Q(\Gamma+\bar{\Gamma})+(L+\bar{L})^*\mathbb{E}^*\bar{Q}\mathbb{E}(\Gamma+\bar{\Gamma})+(\hat{L}+\hat{\bar{L}})^*G_T(\hat{\Gamma}+\hat{\bar{\Gamma}})+(\hat{L}+\hat{\bar{L}})^*\mathbb{E}^*\bar{G}_T\mathbb{E}(\hat{\Gamma}+\hat{\bar{\Gamma}}), \\ [2mm]
\Theta_3=(\Gamma+\bar{\Gamma})^*Q(\Gamma+\bar{\Gamma})+(\Gamma+\bar{\Gamma})^*\mathbb{E}^*\bar{Q}\mathbb{E}(\Gamma+\bar{\Gamma})+ (\hat{\Gamma}+\hat{\bar{\Gamma}})^*G_T(\hat{\Gamma}+\hat{\bar{\Gamma}})+(\hat{\Gamma}+\hat{\bar{\Gamma}})^*\mathbb{E}^*\bar{G}_T\mathbb{E}(\hat{\Gamma}+\hat{\bar{\Gamma}}),
\end{array}
\end{eqnarray*}
and the inner products are understood from the context. Hence, for any $\zeta$, $u\mapsto J(\zeta,u)$ is a quadratic functional on the Hilbert space $\mathcal{U}_{ad}$, and the original problem (MF-LQ) is transformed to a minimization problem of a functional  over $\mathcal{U}_{ad}$. We then have the following results.

\begin{proposition}\label{proposition-operator-Hilbert}
(i). If $J(\zeta, u)$ has a minimum, then
$$\Theta_1\geq 0. $$

(ii). Problem (MF-LQ) is (uniquely) solvable if and only if $\Theta_1\geq 0$ and there  exists a (unique) ${u}$ such that
\begin{eqnarray*}
\Theta_1u+\Theta_2\zeta=0.
\end{eqnarray*}

(iii). If $\Theta_1>0$, then for any $\zeta$, $J(\zeta, u)$ admits a pathwise unique minimizer $u^o$ given by
\begin{eqnarray}\label{proposition-operator-Hilbert-control}
u^o_k=-(\Theta_1^{-1}\Theta_2\zeta)(k), ~~k\in\mathbb{N}.
\end{eqnarray}
In addition, if
\begin{eqnarray}\label{proposition-operator-Hilbert-condition}
Q_k, Q_k+\bar{Q}_k\geq 0,~~R_k, R_k+\bar{R}_k>0,~k\in\mathbb{N},~~G_N,G_N+\bar{G}_N\geq 0,
\end{eqnarray}
then $\Theta_1>0$.
\end{proposition}

Proof. The proofs of \emph{(i)}, \emph{(ii)} and the first part of \emph{(iii)} are well known, and are omitted here; readers may refer to \cite{Yong-MFF-2011}\cite{Yong2012report}\cite{Yong-Zhou} for solutions of similar problems of quadratic functional optimization on Hilbert space. Clearly, from (\ref{proposition-operator-Hilbert-condition}), we have that $(L+\bar{L})^*Q(L+\bar{L})+(L+\bar{L})^*\mathbb{E}^*Q\mathbb{E}(L+\bar{L}) +(\hat{L}+\hat{\bar{L}})^*G(\hat{L}+\hat{\bar{L}})\geq 0$. Also for any nonzero $u\in\mathcal{U}_{ad}$,
\begin{eqnarray*}
\begin{array}{rl}
&\langle Ru, u\rangle+\langle \mathbb{E}^*\bar{R}\mathbb{E}u,u\rangle =\ds\sum_{k=0}^{N-1}\mathbb{E}\left[ u^T_kR_ku_k+(\mathbb{E}u_k)^T\bar{R}_k\mathbb{E}u_k \right] \\ [2mm]
= & \ds\sum_{k=0}^{N-1}\left[\mathbb{E}\left[\left(u_k-\mathbb{E}u_k\right)^TR_k\left(u_k-\mathbb{E}u_k\right) \right]+\left(\mathbb{E}u_k\right)^T(R_k+\bar{R}_k)\left(\mathbb{E}u_k\right)\right]>0.
\end{array}
\end{eqnarray*}
This implies that $R+\mathbb{E}^*\bar{R}\mathbb{E}$ is positive definite. Therefore, $\Theta_1$ is positive definite. This completes the proof.
\hfill $\square$

\begin{remark}Proposition \ref{proposition-operator-Hilbert}
presents necessary and sufficient conditions for the solvability of Problem (MF-LQ); it also shows that the optimal control (\ref{proposition-operator-Hilbert-control}) is a linear functional of initial value $\zeta$. Note that for any $k\in\mathbb{N}$, $u_k$ does not depend on the current state $x_k$ explicitly, and it depends on $\zeta$ and $k$. Therefore, the optimal control (\ref{proposition-operator-Hilbert-control}) may be viewed as an open-loop control.
\end{remark}

We shall show that under condition (\ref{proposition-operator-Hilbert-condition}), the unique optimal control (\ref{proposition-operator-Hilbert-control}) is, indeed, a linear feedback of current state, i.e., a closed-loop control. First, introduce several sequences of operators
\begin{eqnarray}\label{operator-0}
\left\{\begin{array}{l}
\mathcal{A}_kx=A_kx+\bar{A}_k\mathbb{E}x, ~~x\in\mathcal{X}_k ~~k\in \mathbb{N}, \\ [2mm]
\mathcal{B}_ku=B_ku+\bar{B}_k\mathbb{E}u, ~~u\in\mathcal{U}_k,~~k\in \mathbb{N}, \\ [2mm]
\mathcal{C}_kx=C_kx+\bar{C}_k\mathbb{E}x, ~~x\in\mathcal{X}_k,~~k\in \mathbb{N}, \\ [2mm]
\mathcal{D}_ku=D_ku+\bar{D}_k\mathbb{E}u, ~~u\in\mathcal{U}_k,~~k\in \mathbb{N}.
\end{array}\right.
\end{eqnarray}
Clearly, $\mathcal{A}_k, \mathcal{B}_k, \mathcal{C}_k, \mathcal{D}_k, k\in \mathbb{N}$, are all bounded linear operators, defined from $\mathcal{X}_k$ and $\mathcal{U}_k$ to $\mathcal{X}_k$ and $\mathcal{U}_k$, respectively. Consequently, (\ref{system}) may be rewritten as
\begin{eqnarray}\label{system-operator}
x_{k+1}=(\mathcal{A}_kx_k+\mathcal{B}_ku_k)+(\mathcal{C}_kx_k+\mathcal{D}_ku_k)w_k.
\end{eqnarray}
Further, the performance functional may be represented as
\begin{eqnarray*}
\begin{array}{rl}
J(\zeta,u)
= & \ds\sum_{k=0}^{N-1}\left(\langle Q_kx_{k},x_{k}\rangle+\langle\bar{Q}_k\mathbb{E}x_{k},\mathbb{E}x_{k}\rangle+\langle R_ku_k, u_k\rangle+\langle\bar{R}_k\mathbb{E}u_{k}, \mathbb{E}u_{k}\rangle\right) \\ [2mm]
& +\langle{G}_Nx_{N},x_{N}\rangle+\langle\bar{G}_N\mathbb{E}x_{N}, \mathbb{E}x_{N}\rangle \\ [2mm]
= & \ds\sum_{k=0}^{N-1}\left(\langle (Q_k+\mathbb{E}^*\bar{Q}_k\mathbb{E})x_{k},x_{k}\rangle+\langle (R_k+\mathbb{E}^*\bar{R}_k\mathbb{E})u_k, u_k\rangle\right)+\langle({G}_N+ \mathbb{E}^*\bar{G}_N\mathbb{E})x_{N},x_{N}\rangle.
\end{array}
\end{eqnarray*}
Here $\langle Q_kx_{k},x_{k}\rangle$ denotes $\mathbb{E}\left(x_k^TQ_k x_k\right)$, $(Q_k+\mathbb{E}^*\bar{Q}_k\mathbb{E})x_{k}=Q_kx_{k}+ \mathbb{E}^*\bar{Q}_k\mathbb{E}x_{k}$, with similar meanings for related notation.
Define
\begin{eqnarray*}
\left\{\begin{array}{l}
\mathcal{Q}_kx=(Q_k+\mathbb{E}^*\bar{Q}_k\mathbb{E})x,~~x\in \mathcal{X}_k,~~k\in \mathbb{N},\\
\mathcal{R}_ku=(R_k+\mathbb{E}^*\bar{R}_k\mathbb{E})u,~~u\in \mathcal{U}_k,~~k\in \mathbb{N},\\
\mathcal{G}_Nx=({G}_N+ \mathbb{E}^*\bar{G}_N\mathbb{E})x,~~x\in \mathcal{X}_N.
\end{array}\right.
\end{eqnarray*}
Then
\begin{eqnarray}\label{cost-operator}
J(\zeta,u)=\sum_{k=0}^{N-1}\left(\langle \mathcal{Q}_kx_{k},x_{k}\rangle+\langle \mathcal{R}_ku_k, u_k\rangle\right)+\langle\mathcal{G}_Nx_{N},x_{N}\rangle.
\end{eqnarray}
Therefore, Problem (MF-LQ) may be rewritten as the following
\begin{eqnarray}\label{Problem (MF-LQ)-operator}
\left\{\begin{array}{l}
\mbox{minimize } (\ref{cost-operator})\\
\mbox{subject to }u\in\mathcal{U}_{ad},\mbox{ with }(x_\cdot, u_\cdot)\mbox{ satisfying }(\ref{system-operator}).
\end{array}\right.
\end{eqnarray}
Clearly, (\ref{Problem (MF-LQ)-operator}) is an operator stochastic LQ problem in discrete time. Operator LQ problems in the deterministic and continuous-time cases have been studied thoroughly. Generally speaking, an operator LQ problem is a problem of infinite dimensional control theory. For general infinite dimensional control theory, readers may refer to \cite{Bensoussan 1992} and related literature. In this paper, by transforming Problem (MF-LQ) to (\ref{Problem (MF-LQ)-operator}), it is possible to obtain the closed-loop form of the optimal control.

Suppose we have a sequence of self-adjoint bounded linear operators $\{\mathcal{P}_k: \mathcal{X}_k\mapsto \mathcal{X}_k; k\in \bar{\mathbb{N}}\}$, which are determined below. Noting that $\mathcal{X}_{l}\subseteq \mathcal{X}_k$, for $l\leq k$, we assume that
\begin{eqnarray}\label{operator-assumption}
\mathcal{P}_k\mathcal{X}_l\equiv\{\mathcal{P}_kx\,|\,x\in \mathcal{X}_l\}\subseteq \mathcal{X}_l,~~l\leq k.
\end{eqnarray}
This will be proved below. By adding $\langle \mathcal{P}_Nx_N, x_N\rangle-\langle\mathcal{P}_0x_0, x_0\rangle$
to both sides of (\ref{cost-operator}), we have
\begin{eqnarray}\label{cost-operator-2}
\begin{array}{l}
J(\zeta,u)+\langle \mathcal{P}_Nx_N, x_N\rangle-\langle\mathcal{P}_0x_0, x_0\rangle \\ [2mm]
\q =\sum_{k=0}^{N-1}\left(\langle \mathcal{Q}_kx_{k},x_{k}\rangle+\langle \mathcal{R}_ku_k, u_k\rangle+\langle \mathcal{P}_{k+1}x_{k+1}, x_{k+1}\rangle-\langle\mathcal{P}_kx_k, x_k\rangle\right)+\langle\mathcal{G}_Nx_{N},x_{N}\rangle.
\end{array}
\end{eqnarray}
To proceed, we need some calculations.
\begin{eqnarray*}
\begin{array}{l}
\langle \mathcal{P}_{k+1}x_{k+1}, x_{k+1}\rangle \\ [2mm]
\q =\langle \mathcal{P}_{k+1}[(\mathcal{A}_kx_k+\mathcal{B}_ku_k)+(\mathcal{C}_kx_k+\mathcal{D}_ku_k)w_k], [(\mathcal{A}_kx_k+\mathcal{B}_ku_k)+(\mathcal{C}_kx_k+\mathcal{D}_ku_k)w_k]\rangle\\
\q =\mathbb{E}\left[(\mathcal{A}_kx_k+\mathcal{B}_ku_k)^T\left(\mathcal{P}_{k+1}(\mathcal{A}_kx_k+\mathcal{B}_ku_k)\right)\right]+
\mathbb{E}\left[(\mathcal{A}_kx_k+\mathcal{B}_ku_k)^T\left(\mathcal{P}_{k+1}(\mathcal{C}_kx_k+\mathcal{D}_ku_k) w_k\right)\right] \\ [2mm]
\qq +\mathbb{E}\left[(\mathcal{C}_kx_k+\mathcal{D}_ku_k)^T\left(\mathcal{P}_{k+1}(\mathcal{A}_kx_k+\mathcal{B}_ku_k) w_k\right)\right]+\mathbb{E}\left[(\mathcal{C}_kx_k+\mathcal{D}_ku_k)^T\left(\mathcal{P}_{k+1}(\mathcal{C}_kx_k+\mathcal{D}_ku_k)w^2_k\right)\right].
\end{array}
\end{eqnarray*}
Clearly, by (\ref{operator-assumption}), we have
\begin{eqnarray*}
\begin{array}{l}
\mathbb{E}\left[ (\mathcal{C}_kx_k+\mathcal{D}_ku_k)^T\left(\mathcal{P}_{k+1}(\mathcal{A}_kx_k+\mathcal{B}_ku_k) w_k\right)\right] \\ [2mm]
\q = \mathbb{E}\left[\mathbb{E}\left(w_k|\mathcal{F}_{k-1}\right) (\mathcal{C}_kx_k+\mathcal{D}_ku_k)^T\left(\mathcal{P}_{k+1}(\mathcal{A}_kx_k+\mathcal{B}_ku_k)\right) \right]=0, \\ [2mm]
\mathbb{E}\left[(\mathcal{A}_kx_k+\mathcal{B}_ku_k)^T\left(\mathcal{P}_{k+1}(\mathcal{C}_kx_k+\mathcal{D}_ku_k)\right)w_k\right] \\ [2mm]
\q = \mathbb{E}\left[\mathbb{E}\left(w_k|\mathcal{F}_{k-1}\right)(\mathcal{A}_kx_k+\mathcal{B}_ku_k)^T\left(\mathcal{P}_{k+1}(\mathcal{C}_kx_k+\mathcal{D}_ku_k)\right) \right]=0, \\ [2mm]
\mathbb{E}\left[ (\mathcal{C}_kx_k+\mathcal{D}_ku_k)^T\left(\mathcal{P}_{k+1}(\mathcal{C}_kx_k+\mathcal{D}_ku_k) w^2_k\right)\right] \\ [2mm]
\q = \mathbb{E}\left[\mathbb{E}\left(w_k^2|\mathcal{F}_{k-1}\right)  (\mathcal{C}_kx_k+\mathcal{D}_ku_k)^T\left(\mathcal{P}_{k+1}(\mathcal{C}_kx_k+\mathcal{D}_ku_k)\right)\right] \\ [2mm]
\q = \mathbb{E}\left[ (\mathcal{C}_kx_k+\mathcal{D}_ku_k)^T\mathcal{P}_{k+1}(\mathcal{C}_kx_k+\mathcal{D}_ku_k)\right].
\end{array}
\end{eqnarray*}
Therefore, it follows that
\begin{eqnarray}\label{cost-operator-3}
\begin{array}{rl}
\langle \mathcal{P}_{k+1}x_{k+1}, x_{k+1}\rangle
= & \mathbb{E}\left[(\mathcal{A}_kx_k+\mathcal{B}_ku_k)^T\mathcal{P}_{k+1}(\mathcal{A}_kx_k+\mathcal{B}_ku_k) + (\mathcal{C}_kx_k+\mathcal{D}_ku_k)^T\mathcal{P}_{k+1}(\mathcal{C}_kx_k+\mathcal{D}_ku_k)\right] \\ [2mm]
= & \langle \mathcal{A}_k^*\mathcal{P}_{k+1}\mathcal{A}_kx_k,x_k\rangle+\langle \mathcal{C}_k^*\mathcal{P}_{k+1}\mathcal{C}_kx_k,x_k\rangle +\langle \mathcal{A}_k^*\mathcal{P}_{k+1}\mathcal{B}_ku_k, x_k\rangle+\langle \mathcal{B}_k^*\mathcal{P}_{k+1}\mathcal{A}_kx_k, u_k\rangle \\ [2mm]
& +\langle \mathcal{C}_k^*\mathcal{P}_{k+1}\mathcal{D}_ku_k, x_k\rangle+\langle \mathcal{D}_k^*\mathcal{P}_{k+1}\mathcal{C}_kx_k, u_k\rangle+ \langle \mathcal{B}_k^*\mathcal{P}_{k+1}\mathcal{B}_ku_k, u_k\rangle+\langle \mathcal{D}_k^*\mathcal{P}_{k+1}\mathcal{D}_ku_k, u_k\rangle.
\end{array}
\end{eqnarray}
Substituting (\ref{cost-operator-3}) in (\ref{cost-operator-2}), we have
\begin{eqnarray}\label{cost-operator-4}
\begin{array}{rl}
J(x_0,u)= & \ds\sum_{k=0}^{N-1}\left[\langle \left(\mathcal{Q}_k+ \mathcal{A}_k^*\mathcal{P}_{k+1}\mathcal{A}_k+\mathcal{C}_k^*\mathcal{P}_{k+1}\mathcal{C}_k-\mathcal{P}_k\right)x_{k},x_{k}\rangle+2\langle \left(\mathcal{B}_k^*\mathcal{P}_{k+1}\mathcal{A}_k+ \mathcal{D}_k^*\mathcal{P}_{k+1}\mathcal{C}_k \right)x_k, u_k\rangle\right. \\ [2mm]
& \qq \left.+\langle \left(\mathcal{R}_k+ \mathcal{B}_k^*\mathcal{P}_{k+1}\mathcal{B}_k+\mathcal{D}_k^*\mathcal{P}_{k+1}\mathcal{D}_k\right)u_k, u_k\rangle\right]+\langle\left(\mathcal{G}_N-\mathcal{P}_N\right)x_{N},x_{N}\rangle+\langle\mathcal{P}_0x_0, x_0\rangle\\
= & \ds\sum_{k=0}^{N-1}\left[\langle\Theta_{0k} x_k, x_k\rangle+2\langle\Theta_{1k}x_k, u_k\rangle+\langle \Theta_{2k}u_k, u_k\rangle  \right]+\langle\left(\mathcal{G}_N-\mathcal{P}_N\right)x_{N},x_{N}\rangle+\langle\mathcal{P}_0x_0, x_0\rangle,
\end{array}
\end{eqnarray}
where for any $k\in\mathbb{N}$
\begin{eqnarray}\label{theta}
\left\{\begin{array}{l}
\Theta_{0k}=\mathcal{Q}_k+ \mathcal{A}_k^*\mathcal{P}_{k+1}\mathcal{A}_k+\mathcal{C}_k^*\mathcal{P}_{k+1}\mathcal{C}_k-\mathcal{P}_k,\\
\Theta_{1k}=\mathcal{B}_k^*\mathcal{P}_{k+1}\mathcal{A}_k+ \mathcal{D}_k^*\mathcal{P}_{k+1}\mathcal{C}_k,\\
\Theta_{2k}=\mathcal{R}_k+ \mathcal{B}_k^*\mathcal{P}_{k+1}\mathcal{B}_k+\mathcal{D}_k^*\mathcal{P}_{k+1}\mathcal{D}_k.
\end{array}\right.
\end{eqnarray}

Consequently, we have the following result:
\begin{proposition}\label{proposition-operator}
Under the condition that
\begin{eqnarray}\label{proposition-operator-condition}
Q_k, Q_k+\bar{Q}_k\geq 0,~~R_k, R_k+\bar{R}_k>0,~k\in\mathbb{N},~~G_N,G_N+\bar{G}_N\geq 0,
\end{eqnarray}
the unique optimal control for Problem (MF-LQ) is
\begin{eqnarray}\label{proposition-operator-control}
u_k^o=-\left(\mathcal{R}_k+ \mathcal{B}_k^*\mathcal{P}_{k+1}\mathcal{B}_k+\mathcal{D}_k^*\mathcal{P}_{k+1}\mathcal{D}_k \right)^{-1}\left(\mathcal{B}_k^*\mathcal{P}_{k+1}\mathcal{A}_k+ \mathcal{D}_k^*\mathcal{P}_{k+1}\mathcal{C}_k\right)x_k,
\end{eqnarray}
where
\begin{eqnarray}\label{Riccati-equation-operator}
\left\{\begin{array}{rl}
\mathcal{P}_k = & \mathcal{Q}_k+ \mathcal{A}_k^*\mathcal{P}_{k+1}\mathcal{A}_k+\mathcal{C}_k^*\mathcal{P}_{k+1}\mathcal{C}_k \\ [2mm]
& -\left(\mathcal{B}_k^*\mathcal{P}_{k+1}\mathcal{A}_k+ \mathcal{D}_k^*\mathcal{P}_{k+1}\mathcal{C}_k\right)^*\left(\mathcal{R}_k+ \mathcal{B}_k^*\mathcal{P}_{k+1}\mathcal{B}_k+\mathcal{D}_k^*\mathcal{P}_{k+1}\mathcal{D}_k \right)^{-1}\left(\mathcal{B}_k^*\mathcal{P}_{k+1}\mathcal{A}_k+ \mathcal{D}_k^*\mathcal{P}_{k+1}\mathcal{C}_k\right),~k\in\mathbb{N}, \\ [2mm]
\mathcal{P}_N = & \mathcal{G}_N.
\end{array}\right.
\end{eqnarray}
\end{proposition}

\emph{Proof}. By (\ref{cost-operator-4}), if (\ref{Riccati-equation-operator}) is well posed, we have that
\begin{eqnarray*}
\begin{array}{rl}
J(\zeta,u)
= & \ds\sum_{k=0}^{N-1}\left[\langle \Theta_{2k}(u_k+\Theta_{2k}^{-1}\Theta_{1k}x_k), (u_k+\Theta_{2k}^{-1}\Theta_{1k}x_k)\rangle+\langle\Theta_{0k}-\Theta_{1k}^*\Theta_{2k}^{-1}\Theta_{1k}x_k, x_k\rangle \right] \\ [2mm]
&+\langle\left(\mathcal{G}_N-\mathcal{P}_N\right)x_N, x_N\rangle+\langle\mathcal{P}_0\zeta, \zeta \rangle \\ [2mm]
= & \ds\sum_{k=0}^{N-1}\left[\langle \Theta_{2k}(u_k+\Theta_{2k}^{-1}\Theta_{1k}x_k), (u_k+\Theta_{2k}^{-1}\Theta_{1k}x_k)\rangle \right]+\langle\mathcal{P}_0\zeta, \zeta \rangle.
\end{array}
\end{eqnarray*}
Therefore, the optimal control is given by $u^o_k=-\Theta_{2k}^{-1}\Theta_{1k}x_k$, $k\in\mathbb{N}$, which is (\ref{proposition-operator-control}). We shall now prove that (\ref{Riccati-equation-operator}) is well posed. First, for any $k\in\mathbb{N}$ and $u_k\in \mathcal{U}_k$ and $u_k\neq 0$, it is clear that
\begin{eqnarray}\label{proposition-operator-positive}
\begin{array}{rl}
\langle \mathcal{R}_ku_k, u_k \rangle
= & \mathbb{E}\left(u_k^TR_ku_k\right)+(\mathbb{E}u_k)^T\bar{R}_k(\mathbb{E}u_k) \\ [2mm]
= & \mathbb{E}\left[\left(u_k-\mathbb{E}u_k\right)^TR_k\left(u_k-\mathbb{E}u_k\right) \right]+\left(\mathbb{E}u_k\right)^T(R_k+\bar{R}_k)\left(\mathbb{E}u_k\right) \\ [2mm]
\geq & \lambda^{(k)}_1\mathbb{E}|u_k-\mathbb{E}u_k|_m^2+\lambda^{(k)}_2|\mathbb{E}u_k|_m^2 \\ [2mm]
= & \lambda^{(k)}_1\left(\mathbb{E}|u_k|_m^2-|\mathbb{E}u_k|_m^2\right)+\lambda^{(k)}_2|\mathbb{E}u_k|_m^2 \\ [2mm]
\geq & \lambda^{(k)}\left(\mathbb{E}|u_k|_m^2-|\mathbb{E}u_k|_m^2+|\mathbb{E}u_k|_m^2\right) \\ [2mm]
= & \lambda^{(k)}||u_k||_m^2.
\end{array}
\end{eqnarray}
Here $|\cdot|_m$ denotes the norm in $\mathbb{R}^m$; $\lambda^{(k)}_1, \lambda^{(k)}_2$ are the smallest eigenvalues of matrix $R_k$ and $R_k+\bar{R}_k$, respectively, and $\lambda^{(k)}=\min\{\lambda^{(k)}_1, \lambda^{(k)}_2\}$; $||\cdot||_m$ is the norm induced by inner product in $\mathcal{U}_k$. By (\ref{proposition-operator-positive}), we know that $\mathcal{R}_k\geq \lambda^{(k)}I$, where $I$ is the identical operator defined on $\mathcal{U}_k$, for each $k\in\mathbb{N}$. Under assumption (\ref{proposition-operator-condition}), we know that $\lambda^{(k)}>0$. Thus we have that $\mathcal{R}_k$ is positive definite. Furthermore, by known results, we have that
\begin{eqnarray*}
\mathcal{R}^*_k=\left({R}_k+ \mathbb{E}^*\bar{R}_k\mathbb{E} \right)^*={R}^*_k+\left(\mathbb{E}^*\bar{R}_k\mathbb{E}\right)^*={R}^*_k+\mathbb{E}^*\bar{R}^*_k\left(\mathbb{E}^*\right)^*
=R_k+\mathbb{E}^*\bar{R}_k\mathbb{E}=\mathcal{R}_k,
\end{eqnarray*}
which means that $\mathcal{R}_k$ is self-adjoint. Similarly, we have that $\mathcal{G}_N$ is self-adjoint and positive semi-definite. Therefore, it follows
\begin{eqnarray}
\begin{array}{rl}
& \langle\left(\mathcal{R}_{N-1}+ \mathcal{B}_{N-1}^*\mathcal{P}_{N}\mathcal{B}_{N-1}+\mathcal{D}_{N-1}^*\mathcal{P}_{N}\mathcal{D}_{N-1}\right)u_{N-1}, u_{N-1}\rangle \\ [2mm]
= & \langle \mathcal{R}_{N-1}u_{N-1}, u_{N-1}\rangle+\langle \mathcal{G}_{N}\left(\mathcal{B}_{N-1}u_{N-1}\right), \mathcal{B}_{N-1}u_{N-1}\rangle+\langle \mathcal{G}_{N}\left(\mathcal{D}_{N-1}u_{N-1}\right), \mathcal{D}_{N-1}u_{N-1}\rangle \\ [2mm]
\geq &\langle \mathcal{R}_{N-1}u_{N-1}, u_{N-1}\rangle\geq \lambda^{(N-1)}||u_k||_m^2.
\end{array}
\end{eqnarray}
Then, $\Theta_{2(N-1)}=\mathcal{R}_{N-1}+ \mathcal{B}_{N-1}^*\mathcal{P}_{N}\mathcal{B}_{N-1}+\mathcal{D}_{N-1}^*\mathcal{P}_{N}\mathcal{D}_{N-1}$is self-adjoint and positive definite, and invertible. Clearly, $\Theta_{2(N-1)}$ is bounded and linear. By the inverse operator theorem, we have that $\Theta_{2(N-1)}^{-1}$ is bounded and linear. Consequently, (\ref{Riccati-equation-operator}) is well-posed for $k=N-1$. By conditions in the paragraph before (\ref{system-Ex}), we known that for any $l\leq N-1$, $z\in\mathcal{X}_{l}$, $K\in \mathbb{R}^{n\times n}$ and $M\in \mathbb{R}^{m\times n}$, we have that $\mathbb{E}^*K\mathbb{E}z\in \mathcal{X}_{l}$ and $\mathbb{E}^*M\mathbb{E}z\in\mathcal{U}_k$. Therefore,
\begin{eqnarray}\label{proposition-operator-adapted-1}
\left(\mathcal{B}_{N-1}^*\mathcal{P}_{N}\mathcal{A}_{N-1}+ \mathcal{D}_{N-1}^*\mathcal{P}_{N}\mathcal{C}_{N-1}\right)z\in \mathcal{U}_{l},~~z\in \mathcal{X}_{l},~~l\leq N-1,
\end{eqnarray}
and
\begin{eqnarray}\label{proposition-operator-adapted-2}
\left(\mathcal{R}_{N-1}+ \mathcal{B}_{N-1}^*\mathcal{P}_{N}\mathcal{B}_{N-1}+\mathcal{D}_{N-1}^*\mathcal{P}_{N}\mathcal{D}_{N-1} \right)z\in \mathcal{U}_{l},~~z\in\mathcal{U}_{l},~~l\leq N-1.
\end{eqnarray}
From (\ref{proposition-operator-adapted-2}), it follows that
\begin{eqnarray}\label{proposition-operator-adapted-3}
\left(\mathcal{R}_{N-1}+ \mathcal{B}_{N-1}^*\mathcal{P}_{N}\mathcal{B}_{N-1}+\mathcal{D}_{N-1}^*\mathcal{P}_{N}\mathcal{D}_{N-1} \right)^{-1}z\in \mathcal{U}_{l},~~z\in\mathcal{U}_{l},~~l\leq N-1.
\end{eqnarray}
By (\ref{proposition-operator-adapted-1}) and (\ref{proposition-operator-adapted-3}), we have that $u^o_{N-1}$ is $\mathcal{F}_{N-1}$-adapted and thus in $\mathcal{U}_{N-1}$. Further, we have
\begin{eqnarray}\label{proposition-operator-adapted-4}
\mathcal{P}_{N-1}z\in \mathcal{X}_l,~~\mbox{ for any }z\in \mathcal{X}_l,~~ l\leq N-1,
\end{eqnarray}
which is (\ref{operator-assumption}) for $k=N-1$.
Similar to (\ref{cost-operator-4}), we have for any $x_{N-1}\in \mathcal{F}_{N-1}$
\begin{eqnarray*}
\begin{array}{rl}
&\mathbb{E}\left[x_{N-1}^TQ_{N-1}x_{N-1}+\left(\mathbb{E}x_{N-1}\right)^T\bar{Q}_{N-1}\mathbb{E}x_{N-1}+u_{N-1}^TR_{N-1}u_{N-1}+\left(\mathbb{E}u_{N-1}\right)^T\bar{R}_{N-1}\mathbb{E}u_{N-1} \right] \\ [2mm]
&+\mathbb{E}\left[u_{N-1}^TR_{N-1}u_{N-1}+\left(\mathbb{E}u_{N-1}\right)^T\bar{R}_{N-1}\mathbb{E}u_{N-1}\right] \\ [2mm]
= & \langle \Theta_{2(N-1)}(u_{N-1}+\Theta_{2(N-1)}^{-1}\Theta_{1(N-1)}x_{N-1}), (u_{N-1}+\Theta_{2(N-1)}^{-1}\Theta_{1(N-1)}x_{N-1})\rangle+\langle \mathcal{P}_{N-1}x_{N-1}, x_{N-1}\rangle \\ [2mm]
\geq & \langle \mathcal{P}_{N-1}x_{N-1}, x_{N-1}\rangle,
\end{array}
\end{eqnarray*}
and the equality is achieved by $u^o_{N-1}$.
This implies that
\begin{eqnarray*}
\begin{array}{rl}
&\langle \mathcal{P}_{N-1}x_{N-1}, x_{N-1}\rangle \\ [2mm]
= & \mathbb{E}\left[x_{N-1}^TQ_{N-1}x_{N-1}+\left(\mathbb{E}x_{N-1}\right)^T\bar{Q}_{N-1}\mathbb{E}x_{N-1}+u_{N-1}^{oT}R_{N-1}u^o_{N-1}+\left(\mathbb{E}u^o_{N-1}\right)^T\bar{R}_{N-1}\mathbb{E}u^o_{N-1} \right]\geq 0.
\end{array}
\end{eqnarray*}
Thus, $\mathcal{P}_{N-1}$ is positive semi-definite. Clearly, by (\ref{Riccati-equation-operator}), $\mathcal{P}_{N-1}$ is linear and self-adjoint. Now, we are able to prove that $\mathcal{P}_{N-1}$ is bounded. As noted by above, $\Theta_{2(N-1)}^{-1}$ is bounded, so we can easily assert that $\mathcal{P}_{N-1}$ is bounded.  Therefore, we may prove by induction that $\{\mathcal{P}_{k}, k\in \mathbb{N}\}$, are all self-adjoint, positive semi-definite and bounded linear operators, and that the optimal control is given by (\ref{proposition-operator-control}) with $u^o_k\in \mathcal{U}_k$, $k\in \mathbb{N}$. This completes the proof. \hfill $\square$

\section{Solution by Riccati equations}

The results presented in previous section are mathematically pleasing, but they are not in a form which can be implemented, as (\ref{proposition-operator-Hilbert-control}) and the operator Riccati difference equation (\ref{Riccati-equation-operator}) are referenced.
However, Proposition \ref{proposition-operator} provides us with the information that the optimal control of Problem (MF-LQ) is a linear state feedback of operator form. Therefore, it is reasonable to assume that the optimal control $u^o$ takes the form
\begin{eqnarray} \label{optimal-control}
u^o_k=L^o_kx_k+\bar{L}^o_k\mathbb{E}x_k,~~k\in\mathbb{N},
\end{eqnarray}
with $L^o_k, \bar{L}^o_k\in \mathbb{R}^{m\times n}$. To compute the optimal feedback gains $L^o_k, \bar{L}^o_k $, we start from a generic linear feedback control
\begin{eqnarray} \label{control-1}
u_k=L_kx_k+\bar{L}_k\mathbb{E}x_k,~~L_k, \bar{L}_k\in \mathbb{R}^{m\times n},~~k\in\mathbb{N}.
\end{eqnarray}
Under (\ref{control-1}), the closed loop system (\ref{system}) becomes
\begin{eqnarray}\label{system-2}
\left\{\begin{array}{rl}
x_{k+1}= & \left[(A_k+B_kL_k)x_k+[B_k\bar{L}_k+\bar{A}_k+\bar{B}_k(L_k+\bar{L}_k)]\mathbb{E}x_k\right] \\ [2mm]
& +\left[(C_k+D_kL_k)x_k+[D_k\bar{L}_k+\bar{C}_k+\bar{D}_k(L_k+\bar{L}_k)]\mathbb{E}x_k\right]w_k, \\ [2mm]
x_0= & \zeta,
\end{array}\right.
\end{eqnarray}
and the cost functional (\ref{cost-finite}) is
\begin{eqnarray}\label{cost-finite-2}
\begin{array}{rl}
J(\zeta,u)
= & \ds\sum_{k=0}^{N-1}\mathbb{E}\left[x_{k}^TQ_kx_{k}+(\mathbb{E}x_{k})^T\bar{Q}_k\mathbb{E}x_{k}+(L_kx_k+\bar{L}_k\mathbb{E}x_k)^TR_k(L_kx_k+\bar{L}_k\mathbb{E}x_k)\right. \\ [2mm]
& \left.+((L_k+\bar{L}_k)\mathbb{E}x_{k})^T\bar{R}_k(L_k+\bar{L}_k)\mathbb{E}x_{k}\right]+\mathbb{E}\left(x_{N}^T{G}_Nx_{N}\right)+(\mathbb{E}x_{N})^T\bar{G}_N\mathbb{E}x_{N} \\ [2mm]
= & \ds\sum_{k=0}^{N-1}\mathbb{E}\left[x_{k}^T\left(Q_k+L_k^TR_kL_k \right)x_{k}+(\mathbb{E}x_{k})^T\bar{\Phi}_k\mathbb{E}x_{k} \right]+\mathbb{E}\left(x_{N}^T{G}_Nx_{N}\right)+(\mathbb{E}x_{N})^T\bar{G}_N\mathbb{E}x_{N} \\ [2mm]
= & \ds\sum_{k=0}^{N-1}\left\{Tr\left[\left(Q_k+L_k^TR_kL_k \right)\mathbb{E}\left(x_{k}x_{k}^T\right) \right]+Tr\left[\bar{\Phi}_k\left(\mathbb{E}x_{k}(\mathbb{E}x_{k})^T \right)\right]\right\} \\ [2mm]
& +\ds Tr\left[{G}_N\mathbb{E}\left(x_{N}x_{N}^T\right)\right]+Tr\left[\bar{G}_N\left(\mathbb{E}x_{N}(\mathbb{E}x_{N})^T\right)\right],
\end{array}
\end{eqnarray}
where
\begin{eqnarray*}
\bar{\Phi}_k=\bar{Q}_k+(L_k+\bar{L}_k)^T\bar{R}_k(L_k+\bar{L}_k)+L_k^TR_k\bar{L}_k+\bar{L}_k^TR_kL_k+\bar{L}_k^TR_k\bar{L}_k.
\end{eqnarray*}
From the form (\ref{control-1}) of the control, we may view $\{(L_k, \bar{L}_k), k\in \mathbb{N}\}$ as the new control input. Also (\ref{cost-finite-2}) reminds us that $\mathbb{E}x_{k}(\mathbb{E}x_{k})^T$, $\mathbb{E}\left(x_{k}x_{k}^T\right)$
may be considered as the new system states. Write $X_k=\mathbb{E}\left(x_{k}x_{k}^T\right), \bar{X}_k=\mathbb{E}x_{k}(\mathbb{E}x_{k})^T$. Then by (\ref{system-2}), we have
\begin{eqnarray}\label{system-3(1)}
\left\{\begin{array}{rl}
X_{k+1}= & (A_k+B_kL_{k})X_k(A_k+B_kL_{k})^T+(A_k+B_kL_{k})\bar{X}_k\left[\bar{A}_k+B_k\bar{L}_{k}+\bar{B}_k(L_{k}+\bar{L}_{k})\right]^T \\ [2mm]
& +\left[\bar{A}_k+B_k\bar{L}_{k}+\bar{B}_k(L_{k}+\bar{L}_{k})\right]\bar{X}_k(A_k+B_kL_{k})^T \\ [2mm]
& +\left[\bar{A}_k+B_k\bar{L}_{k}+\bar{B}_k(L_{k}+\bar{L}_{k})\right]\bar{X}_k\left[\bar{A}_k+B_k\bar{L}_{k}+\bar{B}_k(L_{k}+\bar{L}_{k})\right]^T \\ [2mm]
& +(C_k+D_kL_{k})X_k(C_k+D_kL_{k})^T+(C_k+D_kL_{k})\bar{X}_k\left[\bar{C}_k+D_k\bar{L}_{k}+\bar{D}_k(L_{l}+\bar{L}_{k})\right]^T \\ [2mm]
& +\left[\bar{C}_k+D_k\bar{L}_{k}+\bar{D}_k(L_{k}+\bar{L}_{k})\right]\bar{X}_k(C_k+D_kL_{k})^T \\ [2mm]
& +\left[\bar{C}_k+D_k\bar{L}_{k}+\bar{D}_k(L_{k}+\bar{L}_{k})\right]\bar{X}_k\left[\bar{C}_k+D_k\bar{L}_{k}+\bar{D}_k(L_{k}+\bar{L}_{k})\right]^T \\ [2mm]
\equiv & \mathcal{X}_k(L_{k},\bar{L}_{k}), \\ [2mm]
X_{0}= & \mathbb{E}(\zeta \zeta^T),
\end{array}\right.
\end{eqnarray}
and
\begin{eqnarray}\label{system-3(2)}
\left\{\begin{array}{l}
\bar{X}_{k+1}=\left[(A_k+\bar{A}_k)+(B_k+\bar{B}_k)(L_{k}+\bar{L}_{k})\right]\bar{X}_k\left[(A_k+\bar{A}_k)+(B_k+\bar{B}_k)(L_{k}+\bar{L}_{k})\right]^T\\
\qq\q\2n\equiv\bar{\mathcal{X}}_k(L_{k},\bar{L}_{k}),\\
\bar{X}_0=\mathbb{E}\zeta \left(\mathbb{E}\zeta\right)^T.
\end{array}\right.
\end{eqnarray}
In the language of $\bar{X}$ and $\bar{X}$, $J(\zeta,u)$ with $u$ defined in (\ref{control-1}) may be represented as%
\begin{eqnarray}
J(\zeta,u)&=&\sum_{k=0}^{N-1}\left\{Tr\left[\left(Q_k+L_k^TR_kL_k \right)X_k \right]+Tr\left(\bar{\Phi}_k\bar{X}_k\right)\right\}+Tr\left({G}_NX_{N}\right)+Tr\left(\bar{G}_N\bar{X}_k\right)\nonumber\\
&\equiv&\mathcal{J}(X_0,\bar{X}_0, \mathcal{L}),
\end{eqnarray}
where $\mathcal{L}\equiv \{L_k, \bar{L}_k, k\in \mathbb{N}\}$. Therefore, Problem (MF-LQ) is equivalent to the following problem:
\begin{eqnarray}\label{problem-matrix-optimization}
\left\{
\begin{array}{l}
\ds\min_{L_k,\bar{L}_k\in R^{m\times n}, k\in \mathbb{N}}\mathcal{J}(X_0,\bar{X}_0, \mathcal{L})\\
\mbox{subject to } (\ref{system-3(1)})(\ref{system-3(2)}).
\end{array}\right.
\end{eqnarray}
Clearly, this is a matrix dynamic optimization problem. A natural
way to deal with this class of problems is by the matrix minimum
principle (\cite{Athans}). Following the framework above, we have
the following results.
\begin{theorem}\label{theorem-finite-horizon}
For {Problem (MF-LQ)}, under the  condition
\begin{eqnarray}\label{proposition-operator-condition-2}
Q_k, Q_k+\bar{Q}_k\geq 0,~~R_k, R_k+\bar{R}_k>0,~k\in\mathbb{N},~~G_N,G_N+\bar{G}_N\geq 0,
\end{eqnarray}
the unique optimal control is
\begin{eqnarray}\label{Theorem-control}
\begin{array}{rl}
u^o_{k} = & -(W_{k}^{(1)})^{-1}H_{k}^{(1)}x_k+\[-(W_{k}^{(2)})^{-1}H_{k}^{(2)}+(W_{k}^{(1)})^{-1}H_{k}^{(1)}\]\mathbb{E}x_k \\ [2mm]
\equiv & L_{k}^ox_k+\bar{L}_{k}^o\mathbb{E}x_k,~~~k\in\mathbb{N}.
\end{array}
\end{eqnarray}
Here,
\begin{eqnarray}\label{Theorem-W-H}
\left\{\begin{array}{l}
W_{k}^{(1)}=R_{k}+B_k^TP_{k+1}B_k+D_k^TP_{k+1}D_k, \\ [2mm]
W_{k}^{(2)}=R_{k}+\bar{R}_{k}+(B_k+\bar{B}_k)^T(P_{k+1}+\bar{P}_{k+1})(B_k+\bar{B}_k)+(D_k+\bar{D}_k)^TP_{k+1}(D_k+\bar{D}_k), \\ [2mm]
H_{k}^{(1)}=B_k^TP_{k+1}A_k+D_k^TP_{k+1}C_k, \\ [2mm]
H_{k}^{(2)}=(B_k+\bar{B}_k)^T(P_{k+1}+\bar{P}_{k+1})(A_k+\bar{A}_k)+(D_k+\bar{D}_k)^TP_{k+1}(C_k+\bar{C}_k),
\end{array}\right.
\end{eqnarray}
with
\begin{eqnarray}\label{Theorem-P}
\left\{\begin{array}{l}
P_{k} = Q_{k}+L_{k}^{oT}R_{k}L^o_{k}+(A_k+B_kL^o_{k})^TP_{k+1}(A_k+B_kL^o_{k}) +(C_k+D_kL^o_{k})^TP_{k+1}(C_k+D_kL^o_{k}),\\
P_N=G_N,
\end{array}
\right.
\end{eqnarray}
\begin{eqnarray}\label{Theorem-(bar-P)}
\left\{
\begin{array}{rcl}
\bar{P}_{k}& = & \bar{Q}_{k}+L_{k}^{oT}R_{k}\bar{L}^o_{k}+\bar{L}_{k}^{oT}R_{k}L^o_{k} +\bar{L}_{k}^{oT}R_{k}\bar{L}^o_{k} \\ [2mm]
&& +(L^o_{k}+\bar{L}^o_{k})^T\bar{R}_{k}(L^o_{k}+\bar{L}^o_{k}) \\ [2mm]
& & +(A_k+B_kL^o_{k})^TP_{k+1}[\bar{A}_k+B_k\bar{L}^o_{k}+\bar{B}_k(L^o_{k}+\bar{L}^o_{k})] \\ [2mm]
&& +[\bar{A}_k+B_k\bar{L}^o_{k}+\bar{B}_k(L^o_{k}+\bar{L}^o_{k})]^TP_{k+1}(A_k+B_kL^o_{k}) \\ [2mm]
&& +[\bar{A}_k+B_k\bar{L}^o_{k}+\bar{B}_k(L^o_{k}+\bar{L}^o_{k})]^TP_{k+1}[\bar{A}_k+B_k\bar{L}^o_{k}+\bar{B}_k(L^o_{k}+\bar{L}^o_{k})] \\ [2mm]
&& +(C_k+D_kL^o_{k})^TP_{k+1}[\bar{C}_k+D_k\bar{L}^o_{k}+\bar{D}_k(L^o_{k}+\bar{L}^o_{k})] \\ [2mm]
&& +[\bar{C}_k+D_k\bar{L}^o_{k}+\bar{D}_k(L^o_{k}+\bar{L}^o_{k})]^TP_{k+1}(C_k+D_kL^o_{k}) \\ [2mm]
&& +[\bar{C}_k+D_k\bar{L}^o_{k}+\bar{D}_k(L^o_{k}+\bar{L}^o_{k})]^TP_{k+1}[\bar{C}_k+D_k\bar{L}^o_{k}+\bar{D}_k(L^o_{k}+\bar{L}^o_{k})] \\ [2mm]
&& +[A_k+\bar{A}_k+(B_k+\bar{B}_k)(L^o_{k}+\bar{L}^o_{k})]^T\bar{P}_{k+1} [A_k+\bar{A}_k+(B_k+\bar{B}_k)(L^o_{k}+\bar{L}^o_{k})],\\
\bar{P}_N&=&\bar{G}_N,
\end{array}
\right.
\end{eqnarray}
having the property
\begin{eqnarray}\label{Theorem-nonnegative}
P_k, P_{k}+\bar{P}_{k}\geq 0,~~k\in\bar{\mathbb{N}}.
\end{eqnarray}
\end{theorem}

\emph{Proof}. Introduce the Lagrangian function associated with
Problem (\ref{problem-matrix-optimization}),
\begin{eqnarray}\label{lagrangian-function}
\begin{array}{rl}
\mathfrak{L} = & \ds\sum_{k=0}^{n-1} \mathfrak{L}_k+Tr\left({G}_NX_{N}\right)+Tr\left(\bar{G}_N\bar{X}_k\right) \\ [2mm]
= & \ds\sum_{k=0}^{n-1} \mathfrak{L}_k+Tr\left[\left(G_N~~\bar{G}_N\right)\left(\begin{array}{c}X_N\\\bar{X}_N \end{array}\right)\right],
\end{array}
\end{eqnarray}
where
\begin{eqnarray}\label{lagrangian-function-2}
\begin{array}{rl}
\mathfrak{L}_k = & Tr\left[\left(Q_k+L_k^TR_kL_k \right)X_k \right]+Tr\left(\bar{\Phi}_k\bar{X}_k\right) \\ [2mm]
& +Tr\left[P_{k+1}\left(\mathcal{X}_k(L_k, \bar{L}_k)-X_{k+1}\right)\right]+Tr\left[\bar{P}_{k+1}\left(\bar{\mathcal{X}}_k(L_k, \bar{L}_k)-\bar{X}_{k+1}\right)\right] \\ [2mm]
= & Tr\left[\left(Q_k+L_k^TR_kL_k \right)X_k \right]+Tr\left(\bar{\Phi}_k\bar{X}_k\right)+Tr\left[\left(P_{k+1}~~\bar{P}_{k+1}\right)\left(\begin{array}{c}\mathcal{X}_k(L_k, \bar{L}_k)-X_{k+1} \\ \bar{\mathcal{X}}_k(L_k, \bar{L}_k)-\bar{X}_{k+1}\end{array}\right) \right],
\end{array}
\end{eqnarray}
and $(P_{k+1}~~ \bar{P}_{k+1}), k\in\bar{\mathbb{N}}$ are the Lagrangian
multipliers. Denote
$\mathbb{{P}}_{k+1}=\left(P_{k+1}~~\bar{P}_{k+1}\right)$,
$\mathbb{X}_k=\left(\begin{array}{c}X_k\\\bar{X}_k
\end{array}\right)$. Clearly, by the matrix minimum
principle (\cite{Athans}), the optimal feedback gains $(L_k^o,
\bar{L}_k^o), k\in\mathbb{N}$ and Lagrangian multipliers
$\mathbb{P}_{k+1}, k\in\mathbb{N}$ satisfy the following first-order
necessary conditions
\begin{eqnarray*}
\left\{
\begin{array}{l}
\frac{\partial\mathfrak{L}_k}{\partial L_k}=0,~~~\frac{\partial\mathfrak{L}_k}{\partial \bar{L}_k}=0,~~\mathbb{P}_k=\frac{\partial\mathfrak{L}_k}{\partial \mathbb{X}_k},~~k\in\mathbb{N},\\
\mathbb{P}_N=\left(G_N~~\bar{G}_N\right),
\end{array}\right.
\end{eqnarray*}
i.e.,
\begin{eqnarray}\label{lagrangian-function-3}
\left\{\begin{array}{l}
\ds\frac{\partial\mathfrak{L}_k}{\partial L_k}=0,~~~\frac{\partial\mathfrak{L}_k}{\partial \bar{L}_k}=0,~~{P}_k=\frac{\partial\mathfrak{L}_k}{\partial {X}_k},~~\bar{P}_k=\frac{\partial\mathfrak{L}_k}{\partial \bar{{X}}_k},~~k\in\mathbb{N}, \\ [2mm]
{P}_N=G_N,~~\bar{P}_N=\bar{G}_N.
\end{array}\right.
\end{eqnarray}
Now, we should calculate several gradient matrices. Noting that for any matrix $Y$, $\frac{\partial}{\partial Y}Tr(AYB)=A^TB^T$ if $AYB$ is meaningful, we have
\begin{eqnarray}\label{lagrangian-function-4}
\begin{array}{rl}
\ds\frac{\partial \mathfrak{L}_k}{\partial L_{k}}
= & \left\{2R_{k}L_{k}X_k+2\left[R_{k}\bar{L}_{k}+\bar{R}_{k}(L_{k}+\bar{L}_{k}) \right]\bar{X}_k\right\} \\ [2mm]
& \ds+2\left\{(B_k+\bar{B}_k)^T\bar{P}_{k+1}(A_k+\bar{A}_k)+(B_k+\bar{B}_k)^T\bar{P}_{k+1}(B_k+\bar{B}_k)(L_{k}+\bar{L}_{k}) \right\}\bar{X}_k \\ [2mm]
& \ds+2\left\{B_k^TP_{k+1}(A_k+B_kL_{k})X_k+B_k^TP_{k+1}(\bar{A}_k+B_k\bar{L}_{k})\bar{X}_k+\bar{B}_k^TP_{k+1}(A_k+B_kL_{k})\bar{X}_k\right. \\ [2mm]
& \ds\left.+B_k^TP_{k+1}\bar{B}_k(L_{k}+\bar{L}_{k})\bar{X}_k+\bar{B}_k^TP_{k+1}(\bar{A}_k+B_k\bar{L}_{k})\bar{X}_k+\bar{B}_k^TP_{k+1}\bar{B}_k(L_{k}+\bar{L}_k)\bar{X}_k\right\} \\ [2mm]
& \ds+2\left\{D_k^TP_{k+1}(C_k+D_kL_{k})X_k+D_k^TP_{k+1}(\bar{C}_k+D_k\bar{L}_{k})\bar{X}_k+\bar{D}_k^TP_{k+1}(C_k+D_kL_{k})\bar{X}_k\right. \\ [2mm]
& \left.+D_k^TP_{k+1}\bar{D}_k(L_{k}+\bar{L}_{k})\bar{X}_k+\bar{D}_k^TP_{k+1}(\bar{C}_k+D_k\bar{L}_{k})\bar{X}_k+\bar{D}_k^TP_{k+1}\bar{D}_k(L_{k}+\bar{L}_{k})\bar{X}_k\right\} \\ [2mm]
\ds = & 2\left[R_kL_{k}+B_k^TP_{k+1}(A_k+B_kL_{k})+D_k^TP_{k+1}(C_k+D_kL_{k}) \right]X_k \\ [2mm]
& \ds+2\left\{\left[R_{k}+\bar{R}_{k}+(B_k+\bar{B}_k)^T(P_{k+1}+\bar{P}_{k+1})(B_k+\bar{B}_k)+(D_k+\bar{D}_k)^TP_{k+1}(D_k+\bar{D}_k)\right]\bar{L}_{k}\right. \\ [2mm]
& \ds+\left[ \bar{R}_{k}+(B_k+\bar{B}_k)^T(P_{k+1}+\bar{P}_{k+1})(B_k+\bar{B}_k)-B_k^TP_{k+1}B_k+\bar{D}_k^TP_{k+1}D_k\right. \\ [2mm]
& \ds\left.+{D}_k^TP_{k+1}\bar{D}_k  +\bar{D}_k^TP_{k+1}\bar{D}_k\right]L_{k}+(B_k+\bar{B}_k)^T(P_{k+1}+\bar{P}_{k+1})(A_k+\bar{A}_k) \\ [2mm]
& \ds\left.+(D_k+\bar{D}_k)^TP_{k+1}(C_k+\bar{C}_k)-B^T_kP_{k+1}A_k-D^T_kP_{k+1}C_k\right\}\bar{X}_k \\ [2mm]
= & \ds2\left(W_k^{(1)}L_k+H_k^{(1)}\right)X_k+2\left(W_{k}^{(2)}\bar{L}_{k}+W_{k}^{3}L_{k}+H_{k}^{(3)}\right)\bar{X}_k \\ [2mm]
= & \ds2\left(W_k^{(1)}L_k+H_k^{(1)}\right)(X_k-\bar{X}_k)+2\left(W_{k}^{(2)}\left(\bar{L}_{k}+L_{k}\right)+H_{k}^{(2)}\right)\bar{X}_k, \\ [2mm]
\end{array}
\end{eqnarray}
\begin{eqnarray}
\begin{array}{rl}
\ds\frac{\partial \mathfrak{L}_k}{\partial \bar{L}_{k}}
= & \ds 2\left\{\left[R_{k}+\bar{R}_{k}+(B_k+\bar{B}_k)^T(P_{k+1}+\bar{P}_{k+1})(B_k+\bar{B}-k)+(D_k+\bar{D}_k)^TP_{k+1}(D_k+\bar{D}_k)\right](\bar{L}_{k}+L_{k})\right. \\ [2mm] \label{lagrangian-function-5}
& \ds\left.+(B_k+\bar{B}_k)^T(P_{k+1}+\bar{P}_{k+1})(A_k+\bar{A}_k)+(D_k+\bar{D}_k)^TP_{k+1}(C_k+\bar{C}_k)\right\}\bar{X}_k \\ [2mm]
= & \ds 2\left(W_{k}^{(2)}\left(L_{k}+\bar{L}_{k}\right)+H_{k}^{(2)}\right)\bar{X}_k.
\end{array}
\end{eqnarray}
Here $W_{k}^{(i)}, H_{k}^{(j)}, i,j=1,2$, are defined in (\ref{Theorem-W-H}), and
\begin{eqnarray*}
\left\{\begin{array}{rl}
W_{k}^{(3)}= & \bar{R}_{k}+(B_k+\bar{B}_k)^T(P_{k+1}+\bar{P}_{k+1})(B_k+\bar{B}_k)-B_k^TP_{k+1}B_k+\bar{D}_k^TP_{k+1}D_k \\ [2mm]
& +{D}_k^TP_{k+1}\bar{D}_k+\bar{D}_k^TP_{k+1}\bar{D}_k, \\ [2mm]
H_{k}^{(3)}= & (B_k+\bar{B}_k)^T(P_{k+1}+\bar{P}_{k+1})(A_k+\bar{A}_k)+(D_k+\bar{D}_k)^TP_{k+1}(C_k+\bar{C}_k) \\ [2mm]
& -B^T_kP_{k+1}A_k-D^T_kP_{k+1}C_k.
\end{array}\right.
\end{eqnarray*}
The following properties are used
\begin{eqnarray*}
W_{k}^{(2)}=W_{k}^{(1)}+W_{k}^{(3)},~~~H_{k}^{(2)}=H_{k}^{(1)}+H_{k}^{(3)}.
\end{eqnarray*}
Combining (\ref{lagrangian-function-3})-(\ref{lagrangian-function-5}), the optimal feedback gains $L^o_{k}$ and $\bar{L}^o_{k}$ must satisfy
\begin{eqnarray}\label{lagrangian-function-7}
\left\{\begin{array}{l}
\left(W_k^{(1)}L_k+H_k^{(1)}\right)(X_k-\bar{X}_k)+\left(W_{k}^{(2)}\left(\bar{L}_{k}+L_{k}\right)+H_{k}^{(2)}\right)\bar{X}_k=0, \\ [2mm]
\left(W_{k}^{(2)}\left(L_{k}+\bar{L}_{k}\right)+H_{k}^{(2)}\right)\bar{X}_k=0.
\end{array}\right.
\end{eqnarray}
Note that (\ref{lagrangian-function-7}) holds for any initial values $X_0-\bar{X}_0=\mathbb{E}\left[(\zeta-\mathbb{E}\zeta)(\zeta-\mathbb{E}\zeta)^T\right]$ and $\bar{X}_0=\mathbb{E}\zeta\left(\mathbb{E}\zeta\right)^T$. Therefore, (\ref{lagrangian-function-7}) reduces to
\begin{eqnarray}\label{lagrangian-function-8}
\left\{\begin{array}{l}
W_{k}^{(1)}L_{k}+H_{k}^{(1)}=0, \\ [2mm]
W_{k}^{(2)}(L_{k}+\bar{L}_{k})+H_{k}^{(2)}=0,
\end{array}\right.
\end{eqnarray}
which is obtained by letting coefficients be zero in (\ref{lagrangian-function-7}).
Clearly, we obtain the optimal feedback gains
\begin{eqnarray*}
\left\{\begin{array}{l}
L^o_{k}=-(W_{k}^{(1)})^{-1}H_{k}^{(1)}, \\ [2mm]
\bar{L}^o_k=-(W_{k}^{(2)})^{-1}H_{k}^{(2)}+(W_{k}^{(1)})^{-1}H_{k}^{(1)}.
\end{array}\right.
\end{eqnarray*}

We now derive the equations that $P_k$ and  $\bar{P}_k$ satisfy. By (\ref{lagrangian-function-3}), we have
\begin{eqnarray*}
\begin{array}{rl}
{P}_k= & \ds\frac{\partial\mathfrak{L}_k}{\partial {X}_k}\Big{|}_{L_k=L_k^o, }=Q_k+(L_k^o)^TR_kL^o_k+(A_k+B_kL^o_k)^TP_{k+1}(A_k+B_kL^o_k) \\ [2mm]
& +(C_k+D_kL^o_k)^TP_{k+1}(C_k+D_kL^o_k), \\ [2mm]
\bar{P}_k= & \ds\frac{\partial\mathfrak{L}_k}{\partial \bar{X}_k}\Big{|}_{L_k=L_k^o, \bar{L}_k=\bar{L}_k^o }=\bar{\Phi}^T_k++(A_k+B_kL^o_{k})^TP_{k+1}\left[\bar{A}_k+B_l\bar{L}^o_{k}+\bar{B}_k(L^o_{k}+\bar{L}^o_{k})\right] \\ [2mm]
& +\left[\bar{A}_k+B_k\bar{L}^o_{k}+\bar{B}_k(L^o_{k}+\bar{L}^o_{k})\right]^T P_{k+1}(A_k+B_kL^o_{k}) \\ [2mm]
& +\left[\bar{A}_k+B_k\bar{L}^o_{k}+\bar{B}_k(L^o_{k}+\bar{L}^o_{k})\right]^T P_{k+1}\left[\bar{A}_k+B_k\bar{L}^o_{k}+\bar{B}_k(L^o_{k}+\bar{L}^o_{k})\right] \\ [2mm]
& +(C_k+D_kL^o_{k})^T P_{k+1}\left[\bar{C}_k+D_k\bar{L}^o_{k}+\bar{D}_k(L^o_{k}+\bar{L}^o_{k})\right] \\ [2mm]
& +\left[\bar{C}_k+D_k\bar{L}^o_{k}+\bar{D}_k(L^o_{k}+\bar{L}^o_{k})\right]^T P_{k+1}(C_k+D_kL^o_{k}) \\ [2mm]
& +\left[\bar{C}_k+D_k\bar{L}^o_{k}+\bar{D}_k(L^o_{k}+\bar{L}^o_{k})\right]^T P_{k+1}\left[\bar{C}_k+D_k\bar{L}^o_{k}+\bar{D}_k(L^o_{k}+\bar{L}^o_{k})\right] \\ [2mm]
& +\left[(A_k+\bar{A}_k)+(B_k+\bar{B}_k)(L^o_{k}+\bar{L}^o_{k})\right]^T\bar{P}_{k+1}\left[(A_k+\bar{A}_k)+(B_k+\bar{B}_k)(L^o_{k}+\bar{L}^o_{k})\right],
\end{array}
\end{eqnarray*}
which are (\ref{Theorem-P}) and (\ref{Theorem-(bar-P)}). The final thing is to assure that for any $k\in\mathbb{N}$, $P_k, P_k+\bar{P}_k\geq 0$. We prove this by induction. Clearly, $P_N, P_N+\bar{P}_N\geq 0$ by definition. For $k=N-1$, $P_{N-1}$ is positive semi-definite, while
\begin{eqnarray*}
\begin{array}{rl}
P_{N-1}+\bar{P}_{N-1} = & Q_{N-1}+\bar{Q}_{N-1}+(L_{N-1}^{o}+\bar{L}_{N-1}^o)^TR_{N-1}(L_{N-1}^{o}+\bar{L}_{N-1}^o) \\ [2mm]
& +[A+\bar{A}+(B+\bar{B})(L^o_{N-1}+\bar{L}^o_{N-1})]^T(P_N+\bar{P}_N) [A+\bar{A}+(B+\bar{B})(L^o_{N-1}+\bar{L}^o_{N-1})] \\ [2mm]
& +[C+\bar{C}+(D+\bar{D})(L^o_{N-1}+\bar{L}^o_{N-1})]^TP_N [C+\bar{C}+(D+\bar{D})(L^o_{N-1}+\bar{L}^o_{N-1})]\geq 0.
\end{array}
\end{eqnarray*}
In addition, by (\ref{Theorem-control}), we have
\begin{eqnarray*}
\begin{array}{rl}
& J_{N-1}^{N}({x}_{N-1}, u^o) \\ [2mm]
= & \mathbb{E}\left(x_{N-1}^TQ_{N-1}x_{N-1}\right)+(\mathbb{E}x_{N-1})^T\bar{Q}_{N-1}\mathbb{E}x_{N-1}+\mathbb{E}\left(u_{N-1}^TR_{N-1}u_{N-1}\right) \\ [2mm]
& +(\mathbb{E}u_{N-1})^T\bar{R}_{N-1}\mathbb{E}u_{N-1}+\mathbb{E}\left(x_{N}^T{G}_Nx_{N}\right)+(\mathbb{E}x_{N})^T\bar{G}_N\mathbb{E}x_{N} \\ [2mm]
= & \mathbb{E}\left[ x_{N-1}^T\left(Q_{N-1}+L_{N-1}^{oT}R_{N-1}L^o_{N-1}\right)x\right] \\ [2mm]
& +(\mathbb{E}x)^T\left[\bar{Q}_{N-1}+L_{N-1}^{oT}R_{N-1}\bar{L}^o_{N-1}+\bar{L}_{N-1}^{oT}R_{N-1}L^o_{N-1} +\bar{L}_{N-1}^{oT}R_{N-1}\bar{L}^o_{N-1}\right. \\ [2mm]
& +\left.(L^o_{N-1}+\bar{L}^o_{N-1})^T\bar{R}_{N-1}(L^o_{N-1}+\bar{L}^o_{N-1})\right]\mathbb{E}x_{N-1} \\ [2mm]
& +\mathbb{E}\[x_{N-1}^T\((A+BL^o_{N-1})^TP_N(A+BL^o_{N-1})+(C+DL^o_{N-1})^TP_N(C+DL^o_{N-1})\)x_{N-1} \] \\ [2mm]
& +(\mathbb{E}x_{N-1})^T\[(A+BL^o_{N-1})^TP_N[\bar{A}+B\bar{L}^o_{N-1}+\bar{B}(L^o_{N-1}+\bar{L}^o_{N-1})] \\ [2mm]
& +[\bar{A}+B\bar{L}^o_{N-1}+\bar{B}(L^o_{N-1}+\bar{L}^o_{N-1})]^TP_N(A+BL^o_{N-1}) \\ [2mm]
& +[\bar{A}+B\bar{L}^o_{N-1}+\bar{B}(L^o_{N-1}+\bar{L}^o_{N-1})]^TP_N[\bar{A}+B\bar{L}^o_{N-1}+\bar{B}(L^o_{N-1}+\bar{L}^o_{N-1})] \\ [2mm]
\end{array}
\end{eqnarray*}
\begin{eqnarray*}
\begin{array}{rl}
& +(C+DL^o_{N-1})^TP_N[\bar{C}+D\bar{L}^o_{N-1}+\bar{D}(L^o_{N-1}+\bar{L}^o_{N-1})] \\ [2mm]
& +[\bar{C}+D\bar{L}^o_{N-1}+\bar{D}(L^o_{N-1}+\bar{L}^o_{N-1})]^TP_N(C+DL^o_{N-1}) \\ [2mm]
& +[\bar{C}+D\bar{L}^o_{N-1}+\bar{D}(L^o_{N-1}+\bar{L}^o_{N-1})]^TP_N[\bar{C}+D\bar{L}^o_{N-1}+\bar{D}(L^o_{N-1}+\bar{L}^o_{N-1})] \\ [2mm]
& +[A+\bar{A}+(B+\bar{B})(L^o_{N-1}+\bar{L}^o_{N-1})]^T\bar{P}_N [A+\bar{A}+(B+\bar{B})(L^o_{N-1}+\bar{L}^o_{N-1})] \]\mathbb{E}x_{N-1} \\ [2mm]
= & \mathbb{E}\left(x_{N-1}^TP_{N-1}x_{N-1}\right)+\left(\mathbb{E}x_{N-1}\right)^T\bar{P}_{N-1}\left(\mathbb{E}x_{N-1}\right).
\end{array}
\end{eqnarray*}
Noting that
\begin{eqnarray*}
J_{k}^N(x_k, u^o|_{\{k, k+1,\cdots,N-1\}})&=&\mathbb{E}\left[\left(x_{k}^TQ_kx_{l}+(\mathbb{E}x_{k})^T\bar{Q}_k\mathbb{E}x_{k}+(u_k^o)TR_ku^o_k+(\mathbb{E}u_{k}^o)^T\bar{R}_k\mathbb{E}u^o_{k}\right) \right]\\
&&+J_{k+1}^N(x_{k+1}, u^o|_{\{k+1, k+2,\cdots,N-1\}}),
\end{eqnarray*}
by induction, we have that $P_k, P_k+\bar{P}_k\geq 0$, for any $k\in\mathbb{N}$. This completes the proof. \hfill $\square$

We may view Problem (MF-LQ) another way. For any $z\in\mathcal{H}=\mathcal{X}_k, \mathcal{U}_k$, we easily see that $\mathbb{E}z$ is orthogonal to $z-\mathbb{E}z$, as
$$\langle \mathbb{E}z, z-\mathbb{E}z\rangle=\mathbb{E}\left[\left(\mathbb{E}z\right)^T(z-\mathbb{E}z)\right]=0.$$
Thus, (\ref{control-1}) is equivalent to
\begin{eqnarray}\label{control-2}
u_k=(L_k+\bar{L}_k)\mathbb{E}x_k+L_k(x_k-\mathbb{E}x_k)\equiv M_k\mathbb{E}x_k+L_k(x_k-\mathbb{E}x_k),
\end{eqnarray}
which is a function of $\mathbb{E}x_k$ and $x_k-\mathbb{E}x_k$. Clearly, in (\ref{control-2}), we may design $M_k$ and $L_k$ independently. This reminds us that we may study Problem (MF-LQ) using coordinates $\{\mathbb{E}x_k, x_k-\mathbb{E}x_k\}$ for any $k\in\bar{\mathbb{N}}$. The system equations that $\mathbb{E}x_k$, $x_k-\mathbb{E}x_k$ satisfy are
\begin{eqnarray}\label{system-Ex-2}
&&\left\{\begin{array}{l}
\mathbb{E}x_{k+1}=(A_k+\bar{A}_k)\mathbb{E}x_k+(B_k+\bar{B}_k)\mathbb{E}u_k,\\
\mathbb{E}x_{0}=\mathbb{E}\zeta,
\end{array}\right. \\
\label{system-x-Ex-2}
&& \left\{\begin{array}{l}
x_{k+1}-\mathbb{E}x_{k+1}=[A_k(x_k-\mathbb{E}x_k)+B_k(u_k-\mathbb{E}u_k)]+[C_k(x_k-\mathbb{E}x_k)+(C_k+\bar{C}_k)\mathbb{E}x_k\\
\qq\qq\qq\q +D_k(u_k-\mathbb{E}u_k)+(D_k+\bar{D}_k)\mathbb{E}u_k]w_k,\\
x_{0}-\mathbb{E}x_{0}=\zeta-\mathbb{E}\zeta.
\end{array}
\right.
\end{eqnarray}
The cost functional $J(\zeta, u)$ may be represented as
\begin{eqnarray}\label{cost-u-Eu}
\begin{array}{rl}
J(\zeta, u)=
& \ds\sum_{k=0}^{N-1}\mathbb{E}\left[\left(x_k-\mathbb{E}x_k\right)^TQ_k\left(x_k-\mathbb{E}x_k\right)+\left(\mathbb{E}x_k\right)^T(Q_k+\bar{Q}_k)\mathbb{E}x_k+\left(u_k-\mathbb{E}u_k\right)^TR_k\left(u_k-\mathbb{E}u_k\right)\right. \\ [5mm]
&\left.+\left(\mathbb{E}u_k\right)^T(R_k+\bar{R}_k)\mathbb{E}u_k\right]+\mathbb{E}\left[\left(x_N-\mathbb{E}x_N\right)^TG_N\left(x_N-\mathbb{E}x_N\right)+\left(\mathbb{E}x_N\right)^T(G_N+\bar{G}_N)\mathbb{E}x_N\right].
\end{array}
\end{eqnarray}

To follow the classical method of completing the square, we introduce two sequences of symmetric matrices $\{S_k, k\in\bar{\mathbb{N}}\}$ and $\{T_k, k\in\bar{\mathbb{N}}\}$, which are determined below. Then,
\begin{eqnarray*}
\begin{array}{rl}
& J(\zeta, u)+\left(x_N-\mathbb{E}x_N\right)^TS_N\left(x_N-\mathbb{E}x_N\right)-\left(x_0-\mathbb{E}x_0\right)^TS_0\left(x_0-\mathbb{E}x_0\right)+\left(\mathbb{E}x_N\right)^TT_N\mathbb{E}x_N-\left(\mathbb{E}x_0\right)^TT_0\mathbb{E}x_0\\
= & \ds\sum_{k=0}^{N-1}\mathbb{E}\left[\left(x_k-\mathbb{E}x_k\right)^TQ_k\left(x_k-\mathbb{E}x_k\right)+\left(\mathbb{E}x_k\right)^T(Q_k+\bar{Q}_k)\mathbb{E}x_k+\left(u_k-\mathbb{E}u_k\right)^TR_k\left(u_k-\mathbb{E}u_k\right)\right.\nonumber \\
& +\left(\mathbb{E}u_k\right)^T(R_k+\bar{R}_k)\mathbb{E}u_k+\left(x_{k+1}-\mathbb{E}x_{k+1}\right)^TS_{k+1}\left(x_{k+1}-\mathbb{E}x_{k+1}\right)-\left(x_k-\mathbb{E}x_k\right)^TS_k\left(x_k-\mathbb{E}x_k\right)\\
& \left.+\left(\mathbb{E}x_{k+1}\right)^T T_{k+1}\mathbb{E}x_{k+1}-\left(\mathbb{E}x_k\right)^TT_k\mathbb{E}x_k\right]+\mathbb{E}\left[\left(x_N-\mathbb{E}x_N\right)^TG_N\left(x_N-\mathbb{E}x_N\right)\right.\\
& \left.+\left(\mathbb{E}x_N\right)^T(G_N+\bar{G}_N)\mathbb{E}x_N\right]\\
= & \ds\sum_{k=0}^{N-1}\mathbb{E}\left[\left(x_k-\mathbb{E}x_k\right)^T\left(Q_k+A_k^TS_{k+1}A_k+C_k^TS_{k+1}C_k\right)\left(x_k-\mathbb{E}x_k\right) \right.\\
& +2\left(x_k-\mathbb{E}x_k\right)^T\left(A_k^TS_{k+1}B_k+C_k^TS_{k+1}D_k-S_k\right)\left(u_k-\mathbb{E}u_k\right)\\
& +\left(u_k-\mathbb{E}u_k\right)^T\left(R_k+B_k^TS_{k+1}B_k+D_k^TS_{k+1}D_k\right)\left(u_k-\mathbb{E}u_k\right)\\
& +\left(\mathbb{E}x_k\right)^T\left(Q_k+\bar{Q}_k+(C_k+\bar{C}_k)^TS_{k+1}(C_k+\bar{C}_k)+(A_k+\bar{A}_k)^TT_{k+1}(A_k+\bar{A}_k)-T_k\right)\mathbb{E}x_k\\
& +2\left(\mathbb{E}x_k\right)^T\left((C_k+\bar{C}_k)^TS_{k+1}(D_k+\bar{D}_k)+(A_k+\bar{A}_k)T_{k+1}(B_k+\bar{B}_k) \right)\mathbb{E}u_k\\
& \left.+\left(\mathbb{E}u_k\right)^T\left(R_k+\bar{R}_k+(D_k+\bar{D}_k)^TS_{k+1}(D_k+\bar{D}_k)+(B_k+\bar{B}_k)^TT_{k+1}(B_k+\bar{B}_k)  \right)\mathbb{E}u_k\right]\\
& +\mathbb{E}\left[\left(x_N-\mathbb{E}x_N\right)^TG_N\left(x_N-\mathbb{E}x_N\right)+\left(\mathbb{E}x_N\right)^T(G_N+\bar{G}_N)\mathbb{E}x_N\right]\\
= & \ds\sum_{k=0}^{N-1}\mathbb{E}\left[\left(x_k-\mathbb{E}x_k\right)^T\left(Q_k+A_k^TS_{k+1}A_k+C_k^TS_{k+1}C_k-S_k-(\bar{H}_k^{(1)})^T(\bar{W}_k^{(1)})^{-1}\bar{H}_k^{(1)}\right)\left(x_k-\mathbb{E}x_k\right)\right.\\
& +\left(u_k-\mathbb{E}u_k+(\bar{W}_k^{(1)})^{-1}\bar{H}_k^{(1)}\left(x_k-\mathbb{E}x_k\right)\right)^T\bar{W}_k^{(1)}\left(u_k-\mathbb{E}u_k+(\bar{W}_k^{(1)})^{-1}\bar{H}_k^{(1)}\left(x_k-\mathbb{E}x_k\right)\right)\\
& +\left(\mathbb{E}x_k\right)^T\left(Q_k+\bar{Q}_k+(C_k+\bar{C}_k)^TS_{k+1}(C_k+\bar{C}_k)+(A_k+\bar{A}_k)^TT_{k+1}(A_k+\bar{A}_k)-T_k\right)\mathbb{E}x_k\\
& -\left(\mathbb{E}x_k\right)^T\left((\bar{H}_k^{(2)})^T(\bar{W}_k^{(2)})^{-1}\bar{H}_k^{(2)}\right)\mathbb{E}x_k\\
& +\left.\left(\mathbb{E}u_k+(\bar{W}_k^{(2)})^{-1}\bar{H}_k^{(2)}\mathbb{E}x_k\right)^T\bar{W}_k^{(2)}\left(\mathbb{E}u_k+(\bar{W}_k^{(2)})^{-1}\bar{H}_k^{(2)}\mathbb{E}x_k\right)\right]\\
& +\mathbb{E}\left[\left(x_N-\mathbb{E}x_N\right)^TG_N\left(x_N-\mathbb{E}x_N\right)+\left(\mathbb{E}x_N\right)^T(G_N+\bar{G}_N)\mathbb{E}x_N\right].
\end{array}
\end{eqnarray*}
Here
\begin{eqnarray}\label{Theorem-W-H-2}
\left\{
\begin{array}{l}
\bar{W}_{k}^{(1)}=R_{k}+B_k^TS_{k+1}B_k+D_k^TS_{k+1}D_k,\\
\bar{H}_{k}^{(1)}=B_k^TS_{k+1}A_k+D_k^TS_{k+1}C_k,\\
\bar{W}_{k}^{(2)}=R_{k}+\bar{R}_{k}+(B_k+\bar{B}_k)^TT_{k+1}(B_k+\bar{B}_k)+(D_k+\bar{D}_k)^TS_{k+1}(D_k+\bar{D}_k),\\
\bar{H}_{k}^{(2)}=(B_k+\bar{B}_k)^TT_{k+1}(A_k+\bar{A}_k)+(D_k+\bar{D}_k)^TS_{k+1}(C_k+\bar{C}_k).
\end{array}
\right.
\end{eqnarray}
Letting
\begin{eqnarray}\label{Riccati-equation-2}
\left\{\begin{array}{l}
S_k=Q_k+A_k^TS_{k+1}A_k+C_k^TS_{k+1}C_k-(\bar{H}_k^{(1)})^T(\bar{W}_k^{(1)})^{-1}\bar{H}_k^{(1)},\\
T_k=Q_k+\bar{Q}_k+(C_k+\bar{C}_k)^TS_{k+1}(C_k+\bar{C}_k)+(A_k+\bar{A}_k)^TT_{k+1}(A_k+\bar{A}_k)
-(\bar{H}_k^{(2)})^T(\bar{W}_k^{(2)})^{-1}\bar{H}_k^{(2)},\\
S_N=G_N,~~T_N=G_N+\bar{G}_N.
\end{array}\right.
\end{eqnarray}
and
\begin{eqnarray*}
\left\{\begin{array}{l}
\bar{u}_k-\mathbb{E}\bar{u}_k=-(\bar{W}_k^{(1)})^{-1}\bar{H}_k^{(1)}\left(x_k-\mathbb{E}x_k\right),~~k\in\mathbb{N},\\
\mathbb{E}\bar{u}_k=-(\bar{W}_k^{(2)})^{-1}\bar{H}_k^{(2)}\mathbb{E}x_k,~~k\in\mathbb{N},
\end{array}\right.
\end{eqnarray*}
that is,
\begin{eqnarray}\label{control-3}
\bar{u}_k=-(\bar{W}_k^{(2)})^{-1}\bar{H}_k^{(2)}\mathbb{E}x_k-(\bar{W}_k^{(1)})^{-1}\bar{H}_k^{(1)}\left(x_k-\mathbb{E}x_k\right),~~k\in\mathbb{N},
\end{eqnarray}
we have
\begin{eqnarray*}
J(\zeta,\bar{u})=\mathbb{E}\left[\left(\zeta-\mathbb{E}\zeta\right)^TS_0\left(\zeta-\mathbb{E}\zeta\right)+\left(\mathbb{E}\zeta\right)^TT_0\mathbb{E}\zeta \right]\leq J(\zeta, u)
\end{eqnarray*}
for any $u=(u_0,\cdots,u_{N-1})\mbox{ with }u_k\in\mathcal{U}_k$.
This means that (\ref{control-3}) is an optimal control. In the
following, we show that (\ref{control-3}) is equal to
(\ref{Theorem-control}). We need firstly to show that for any $k\in
\bar{\mathbb{N}}$, $P_k=S_k, T_k=P_k+\bar{P}_k$. In fact, by
substituting (\ref{Theorem-control}) into (\ref{Theorem-P}), we have
\begin{eqnarray}\label{Theorem-P-2}
\begin{array}{rl}
P_k=& Q_k+A_k^TP_{k+1}A_k+C^T_kP_{k+1}C_k+\left(L_k^o\right)^T(R_k+B_k^TP_{k+1}B_k+D_k^TP_{k+1}D_k)L_k^o\\
& +\left(L_k^o\right)^TH_k^{(1)}+\left(H_k^{(1)}\right)^TL_k^o \\
=& Q_k+A_k^TP_{k+1}A_k+C^T_kP_{k+1}C_k-({H}_k^{(1)})^T({W}_k^{(1)})^{-1}{H}_k^{(1)}.
\end{array}
\end{eqnarray}
Noting that $P_N=G_N=S_N$, by (\ref{Theorem-W-H})(\ref{Theorem-W-H-2})(\ref{Riccati-equation-2})(\ref{Theorem-P-2}), we have that
\begin{eqnarray*}
S_k=P_k,~~\bar{W}_{k}^{(1)}={W}_{k}^{(1)},~~\bar{H}_{k}^{(1)}={H}_{k}^{(1)},~~k\in\bar{\mathbb{N}}.
\end{eqnarray*}
Similarly, we have by (\ref{Theorem-(bar-P)})(\ref{Theorem-P})
\begin{eqnarray*}
\begin{array}{rl}
\bar{P}_k = & \bar{Q}_k-L_k^{oT}R_kL_k^o-\left(A_k+B_kL_k^o\right)^TS_{k+1}\left(A_k+B_kL_k^o\right)-\left(C_k+D_kL_k^o\right)^TS_{k+1}\left(C_k+D_kL_k^o\right)\\
&+(C_k+\bar{C}_k)^TS_{k+1}(C_k+\bar{C}_k)+(A_k+\bar{A}_k)^T(P_{k+1}+\bar{P}_{k+1})(A_k+\bar{A}_k)-({H}_k^{(2)})^T({W}_k^{(2)})^{-1}{H}_k^{(2)}\\
= & Q_k+\bar{Q}_k-P_k+(C_k+\bar{C}_k)^TS_{k+1}(C_k+\bar{C}_k)+(A_k+\bar{A}_k)^T(P_{k+1}+\bar{P}_{k+1})(A_k+\bar{A}_k)\\
&-({H}_k^{(2)})^T({W}_k^{(2)})^{-1}{H}_k^{(2)}.
\end{array}
\end{eqnarray*}
Therefore,
\begin{eqnarray}
P_k+\bar{P}_k&=&Q_k+\bar{Q}_k+(C_k+\bar{C}_k)^TS_{k+1}(C_k+\bar{C}_k)+(A_k+\bar{A}_k)^T(P_{k+1}+\bar{P}_{k+1})(A_k+\bar{A}_k)\nonumber\\
&&-({H}_k^{(2)})^T({W}_k^{(2)})^{-1}{H}_k^{(2)}.
\end{eqnarray}
Comparing this with (\ref{Riccati-equation-2}), as $T_N=G_N+\bar{G}_N=P_N+\bar{P}_N$, we have that
\begin{eqnarray*}
T_k=P_k+\bar{P}_k,~~ \bar{W}_{k}^{(2)}={W}_{k}^{(2)},~~\bar{H}_{k}^{(2)}={H}_{k}^{(2)},~~k\in \bar{\mathbb{N}}.
\end{eqnarray*}
Therefore,  (\ref{control-3})  equals to (\ref{Theorem-control}).

To summarize, Theorem \ref{theorem-finite-horizon} may be rewritten in a more compact form as the following results:

\begin{theorem}\label{theorem-finite-horizon2}
For {Problem (MF-LQ)}, under the condition
\begin{eqnarray}\label{proposition-operator-condition-3}
Q_k, Q_k+\bar{Q}_k\geq 0,~~R_k, R_k+\bar{R}_k>0,~k\in\mathbb{N},~~G_N,G_N+\bar{G}_N\geq 0,
\end{eqnarray}
the optimal control is
\begin{eqnarray}\label{Theorem-control2}
{u}^o_k&=&-({W}_k^{(2)})^{-1}{H}_k^{(2)}\mathbb{E}x_k-({W}_k^{(1)})^{-1}{H}_k^{(1)}\left(x_k-\mathbb{E}x_k\right)\equiv M_{k}^o\mathbb{E}x_k+{L}_{k}^o(x_k-\mathbb{E}x_k).
\end{eqnarray}
Here,
\begin{eqnarray}\label{Theorem-W-H2}
\left\{\begin{array}{l}
W_{k}^{(1)}=R_{k}+B_k^TS_{k+1}B_k+D_k^TS_{k+1}D_k,\\
W_{k}^{(2)}=R_{k}+\bar{R}_{k}+(B_k+\bar{B}_k)^TT_{k+1}(B_k+\bar{B}_k)+(D_k+\bar{D}_k)^TS_{k+1}(D_k+\bar{D}_k),\\
H_{k}^{(1)}=B_k^TS_{k+1}A_k+D_k^TS_{k+1}C_k,\\
H_{k}^{(2)}=(B_k+\bar{B}_k)^TT_{k+1}(A_k+\bar{A}_k)+(D_k+\bar{D}_k)^TS_{k+1}(C_k+\bar{C}_k),\\
S_k=Q_k+A_k^TS_{k+1}A_k+C_k^TS_{k+1}C_k-({H}_k^{(1)})^T({W}_k^{(1)})^{-1}{H}_k^{(1)},\\
T_k=Q_k+\bar{Q}_k+(C_k+\bar{C}_k)^TS_{k+1}(C_k+\bar{C}_k)+(A_k+\bar{A}_k)^TT_{k+1}(A_k+\bar{A}_k)
-({H}_k^{(2)})^T({W}_k^{(2)})^{-1}{H}_k^{(2)},\\
S_N=G_N,~~T_N=G_N+\bar{G}_N.
\end{array}\right.
\end{eqnarray}
with the following property
\begin{eqnarray}\label{Theorem-nonnegative2}
P_k, T_k\geq 0,~~k\in\bar{\mathbb{N}}.
\end{eqnarray}
\end{theorem}

\section{Numerical results}

We consider a $4$-period numerical example
$$
\begin{array}{rl}
\ds\min_{u_0,u_1,u_2,u_3} & \ds\mathbb{E}\left[\sum_{k=0}^3\left(x_{k}^TQ_kx_{k}+(\mathbb{E}x_{k})^T\bar{Q}_k\mathbb{E}x_{k}+u_k^TR_ku_k+(\mathbb{E}u_{k})^T\bar{R}_k\mathbb{E}u_{k}\right) \right] \\ [5mm]
& \quad +\mathbb{E}\left(x_4^T{G}_4x_4\right)+(\mathbb{E}x_4)^T\bar{G}_4\mathbb{E}x_4, \\ [5mm]
\mbox{subject to} & \left\{\begin{array}{ll}
x_{k+1}=(A_kx_k+\bar{A}_k\mathbb{E}x_k+B_ku_k+\bar{B}_k\mathbb{E}u_k)+(C_kx_k+\bar{C}_k\mathbb{E}x_k+D_ku_k+\bar{D}_k\mathbb{E}u_k)w_k, \\
x_0 \in \mathbb{R}^3, \hfill k = 0,1,2,3,
\end{array}\right.\end{array}
$$
with coefficients for $k = 0,1,2,3$ as follows
$$\begin{array}{lll}
A_k = \left[\begin{array}{crr}
0.2   &  0.4   &  0.2  \\
         0  &  0.2   &  0.6  \\
    0.6   &  0.4   &  0.2
\end{array}\right],
& \quad\quad &
\bar A_k = \left[\begin{array}{crr}
    0.3&    0.4&    0.2 \\
    0  &    0.2&    0.7 \\
    0.6&    0.5&    0.2
\end{array}\right], \\ [5mm]
B_k = \left[\begin{array}{rr}
    0.4&    0.2 \\
    0.2&    0.4 \\
    0.3&    0.3
\end{array}\right],
& &
\bar B_k = \left[\begin{array}{rr}
    0.5&    0.2 \\
    0.2&    0.5 \\
    0.2&    0.3
\end{array}\right],
\end{array}
$$
$$\begin{array}{lll}
C_k = \left[\begin{array}{rrr}
    0.2&    0.4&    0.6 \\
    0.4&    0.2&    0.6 \\
    0.2&    0.4&    0.2
\end{array}\right],
& &
\bar C_k = \left[\begin{array}{rrr}

    0.3&    0.4&    0.6 \\
    0.4&    0.3&    0.6 \\
    0.2&    0.4&    0.3
\end{array}\right], \\ [5mm]
D_k = \left[\begin{array}{rr}
    0.2&    0.6 \\
    0.6&    0.4 \\
    0.3&    0.1
\end{array}\right],
& &
\bar D_k = \left[\begin{array}{rr}
    0.3&    0.5 \\
    0.5&    0.4 \\
    0.3&    0.3
\end{array}\right], \\ [5mm]
Q_k = \mbox{diag}([0,1.5,1]), & & \bar Q_k = \mbox{diag}([1,1,0]), \\
R_k = \mbox{diag}([1,1]), & & \bar R_k = \mbox{diag}([1.5,1]), \\
G_4 = \mbox{diag}([0,1,1]), & & \bar G_4 = \mbox{diag}([0.5,1,0]).
\end{array}
$$
Based on the Riccati equations (\ref{Theorem-W-H2}), we have
Riccati solutions for $S_i$ and $T_i$ for $i = 0,1,2,3$ are given by
$$
\begin{array}{lll}
S_0 = \left[\begin{array}{rrr}
0.5227  &  0.3542 &   0.1966 \\
    0.3542  &  1.9655  &  0.3170 \\
    0.1966  &  0.3170  &  1.7009
\end{array}\right],
& \quad\quad &
T_0 = \left[\begin{array}{rrr}
4.3329  &  1.7927 &  -0.2507 \\
    1.7927  &  4.4213  &  0.4463 \\
   -0.2507  &  0.4463  &  3.4720
\end{array}\right], \\ [5mm]
S_1= \left[\begin{array}{rrr}
0.5188  &  0.3513  &  0.1951 \\
    0.3513  &  1.9595  &  0.3130 \\
    0.1951  &  0.3130  &  1.6943
\end{array}\right],
& &
T_1 = \left[\begin{array}{rrr}
4.2341  &  1.7868  & -0.2366 \\
    1.7868  &  4.4007 &   0.4611 \\
   -0.2366  &  0.4611 &   3.4411
\end{array}\right], \\ [5mm]
S_2 = \left[\begin{array}{rrr}
0.4862  &  0.3264  &  0.1861 \\
    0.3264  &  1.9219  &  0.2928 \\
    0.1861 &  0.2928  &  1.6660
\end{array}\right],
& &
T_2 = \left[\begin{array}{rrr}
3.4908  &  1.6119  & -0.0389 \\
    1.6119  &  4.2394  &  0.4881 \\
   -0.0389  &  0.4881  &  3.3283
\end{array}\right], \\ [5mm]
S_3 = \left[\begin{array}{rrr}
0.3747  &  0.2421  &  0.1492 \\
    0.2421 &   1.7652  &  0.1849 \\
    0.1492 &   0.1849  &  1.4532 \\
\end{array}\right],
& &
T_3 = \left[\begin{array}{rrr}
1.4782  &  0.3777  &  0.3734 \\
    0.3777  &  3.2001  &  0.5548 \\
    0.3734  &  0.5548  &  2.4932
\end{array}\right].
\end{array}
$$
Then applying Theorem \ref{theorem-finite-horizon2}, we get the
optimal control below
$$
u_k^o = M_{k}^o\mathbb{E}x_k+{L}_{k}^o(x_k-\mathbb{E}x_k), \quad k = 0,1,2,3,
$$
where
$$
\begin{array}{lll}
M_0^o = \left[\begin{array}{rrr}
-0.3286 &  -0.4234  & -0.3474 \\
-0.3189 &  -0.4351  & -0.7770
\end{array}\right],
& \quad\quad &
L_0^o = \left[\begin{array}{rrr}
-0.3455 &  -0.3271  & -0.4240 \\
-0.2467 &  -0.2937  & -0.4941
\end{array}\right], \\ [5mm]
M_1^o = \left[\begin{array}{rrr}
-0.3436 &  -0.4156  & -0.3531 \\
-0.3137 &  -0.4381  & -0.7687
\end{array}\right],
& &
L_1^o = \left[\begin{array}{rrr}
-0.3436 &  -0.3235  & -0.4207 \\
-0.2446 &  -0.2897  & -0.4885
\end{array}\right], \\ [5mm]
M_2^o = \left[\begin{array}{rrr}
-0.4029 &  -0.3946  & -0.3315 \\
-0.2938 &  -0.4160  & -0.7519
\end{array}\right],
& &
 L_2^o = \left[\begin{array}{rrr}
-0.3290 &  -0.3009  & -0.4043 \\
-0.2298 &  -0.2692  & -0.4650
\end{array}\right], \\ [5mm]
M_3^o = \left[\begin{array}{rrr}
-0.2418 &  -0.2552  & -0.3178 \\
-0.1351 &  -0.2213  & -0.5101
\end{array}\right],
& &
 L_3^o = \left[\begin{array}{rrr}
-0.2552 &  -0.2084  & -0.2954 \\
-0.1744 &  -0.1608  & -0.3028
\end{array}\right].
\end{array}
$$

\section{Conclusion}

In this paper, we give four methods to deal with the discrete time mean-field LQ problem: the quadratic optimization method in Hilbert space, the operator LQ method, the matrix dynamic optimization method and the method of completing the square. The optimal control is a linear state feedback using two Riccati equation. For future research, we may consider an infinite horizon mean-field LQ problem. In that case, the stability of the system should first be considered. We shall investigate this in the near future.

%\begin{figure}[h]
% \centering
% \includegraphics[height=0.275\textwidth]{Example-4-1.eps}\\
% \centerline{Fig.7.  $\sum_{i,j=1,2,3,4}|x_t^j-x_t^i|$ of Example 6.4 }\end{figure}
%\vspace{0.5em}

%\begin{figure}[h]
% \centering
% \includegraphics[height=0.275\textwidth]{Example-4-2.eps}\\
%  \centerline{Fig.8.  Curves of the states of Example 6.4 }
%\end{figure}

%\bibliographystyle{plain}        % Include this if you use bibtex
%\bibliography{autosam}           % and a bib file to produce the
                                 % bibliography (preferred). The
                                 % correct style is generated by
                                 % Elsevier at the time of printing.

\end{document}